\documentclass[a4paper,11pt]{amsart}
\usepackage{amssymb,amsfonts,amsmath}
\def\d{\delta}

\def\C{\mathbb{C}}
\def\c2{\mathbb{C}^2}
\def\R{\mathbb{R}}
\def\Q{\mathbb{Q}}
\def\Z{\mathbb{Z}}
\def\N{\mathbb{N}}

\def\1{\bold{1}}

\def\a{\alpha}

\def\e{\varepsilon}
\def\l{\lambda}
\def\f{\varphi}
\def\g{\gamma}

\def\p{\psi}

\def\o{\omega}

\def\Ox{\mathcal{O}_V}
\def\Oxr{\mathcal{O}_{V^{reg}}}

\newtheorem{lem}{Lemma}[section]
\newtheorem{pro}[lem]{Proposition}
\newtheorem{defi}[lem]{Definition}
\newtheorem{def/not}[lem]{Definition/Notations}

\newtheorem{thm}[lem]{Theorem}
\newtheorem{cor}[lem]{Corollary}
\newtheorem{rqe}[lem]{Remark}
\newtheorem{rqes}[lem]{Remarks}
\newtheorem{exa}[lem]{Example}
\newtheorem{exas}[lem]{Examples}

\newenvironment{proof2.1}
{\noindent {\it{Proof of theorem 2.1}}}{$\Box$ \linebreak[4]}

\newenvironment{ack}
{\vskip .2cm \noindent {{\bf Acknowledgements.}}}{\hfill\break}

\begin{document}

\title[Singular K\"ahler-Einstein metrics]
{Singular K\"ahler-Einstein metrics}

\author{Philippe EYSSIDIEUX \& Vincent GUEDJ \& Ahmed ZERIAHI}

\begin{abstract} 
We study degenerate complex Monge-Amp\`ere equations of the form 
$(\omega+dd^c\f)^n = e^{t \f}\mu$
where $\omega$ is a big semi-positive form on a compact K\"ahler  manifold $X$ of 
dimension $n$, $t \in \R^+$, and $\mu=f\omega^n$ is a positive measure with density
$f\in L^p(X,\omega^n)$, $p>1$. We prove the existence and unicity of bounded 
$\o$-plurisubharmonic solutions. We also prove that the solution is continuous under a further technical condition.

In case $X$ is projective and $\omega=\psi^*\omega'$, where $\psi:X\to V$ is a proper
birational morphism to a normal projective variety, $[\omega']\in NS_{\R} (V)$ 
is an ample class  and $\mu$ has only algebraic singularities,
we prove that the solution is smooth in the regular locus of the equation. 

We use these results to construct singular K\"ahler-Einstein metrics of 
non-positive curvature 
on projective klt pairs, in particular on 
canonical models of algebraic varieties of general type. 
\end{abstract}

\maketitle

{ 2000 Mathematics Subject Classification:} {\it 32W20, 32Q20, 32J27, 14J17}.

\section*{Introduction}

Thirty years ago, in a celebrated article [Y],  S.T. Yau (and independently T. Aubin [A]) 
solved the Calabi conjecture by studying complex
Monge-Amp\`ere equations on a compact K\"ahler manifold.

Since then, complex Monge-Amp\`ere equations have been extremely useful in  
K\"ahler geometry (see for instance [DP]) and 
in the dynamical study of rational mappings (see [S] and references therein). 

Two major developments in the theory of complex Monge-Amp\`ere equations occurred in the last 
decade. In the local theory, a deeper analysis of the
image of the complex Monge-Amp\`ere operator [C], [K 1] has followed the 
pioneering work of E.Bedford and A.Taylor [BT].
In the global theory a new proof of the ${\mathcal C}^0$-estimate [K 1,2]
has allowed one to treat complex Monge-Amp\`ere equations with more degenerate R.H.S.

In [GZ1], [GZ2], two of us revisited and extended the results of [BT], [C],
on complex Monge-Amp\`ere operators to compact K\"ahler manifolds. 
In the present article, we use these methods to study complex Monge-Amp\`ere equations 
with degenerate L.H.S. We first define, in the spirit of [C], [GZ 2],
weak solutions to degenerate complex Monge-Amp\`ere equations and then
prove, using ideas of [K 1,2], that these solutions are bounded:
\vskip.2cm

\noindent{\bf{Theorem A.}} 
{\em  Let $X$ be a compact K\"ahler manifold, $\omega$ a semi positive (1,1)-form such that 
$\int_X \omega^n >0$ and
$0 \leq f\in L^p (X,\omega^n)$, $p>1$, a density such that $\int_X f \o^n=\int_X \o^n$.
Then there is a unique bounded function $\f$ on $X$
such that $\o+dd^c \f \geq 0$ and
$$
(\omega+dd^c\f)^n = f \o^n \; \; 
\text{ with } \; \; \sup_X \f=-1.
$$
 
Furthermore, $\f$ is continuous provided there exists a decreasing sequence $\f_j$ of continuous $\omega$-psh functions with $\lim \f_j= \f$ and $f \mapsto \f$ is a continuous map from $L^p(X,\omega^n)$ to $L^{\infty}(X)$.}

\vskip.2cm

When $\omega$ has algebraic singularities, then $\mu$ can be assumed to have 
$L^p$ density ($p>1$) with respect to the Lebesgue measure. 

When for every $\omega$-psh $\f$
there is a decreasing sequence of continuous $\omega$-psh functions converging to $\f$ pointwise, we will say $(X,[\omega])$ enjoys the continuous approximation property. We believe it holds 
in great generality but the problem turned out to be more subtle than expected- hence we leave this as an open problem for further research.

With this $L^{\infty}$-estimate, it is possible to adapt classical ideas of 
[Y] and [Ts] and prove: 
\vskip.2cm

\noindent {\bf{Theorem B.}}
{\em Let $X$ be projective algebraic complex manifold, $\omega$ a 
smooth semi positive closed $(1,1)$-form that is positive outside a complex subvariety $S \subset X$.
Let $\Omega$ be a K\"ahler form on $X$.
Assume furthermore  that $\omega^n= D \Omega^n$ where $D^{-\epsilon}$ 
is in $L^1 (\Omega^n)$ and that $[\omega], [\Omega] \in NS_{\R}(X)$. 

 Let $\sigma_1, ..., \sigma_p$ (resp. $\tau_1, ..., \tau_q$) be holomorphic sections of some line bundle. 
 (resp. of  some other line bundle). Assume $k,l\in \R_{\ge 0}$ 
and $F \in {\mathcal C}^{\infty}(X,\R)$ are fixed so that
$$
\int_X \frac{1}{|\tau_1|^{2l}+ \ldots + |\tau_q|^{2l}}  \Omega^n < \infty  
\text{ and }
\int_X \frac{|\sigma_1|^{2k}+ \ldots + |\sigma_p|^{2k}}{|\tau_1|^{2l}+ \ldots + |\tau_q|^{2l}} 
e^F \Omega^n= \int_X \omega^n.
$$

Then the unique bounded function $\f$ such that $\o+dd^c \f \geq 0$ and
$$
(\omega+ dd^c \f)^n = \frac{|\sigma_1|^{2k}+ \ldots + |\sigma_p|^{2k}}{|\tau_1|^{2l}+ \ldots + |\tau_q|^{2l}} 
e^{F} \Omega^n,
\text{ with } \sup_X \f=-1,
$$
is smooth outside $B=S\cup \cap_i \{\sigma_i=0\}\cup \cap_i \{\tau_i=0\} $. }
\vskip.2cm

 This result should be compared with
[Y], Theorem 8, p. 403. Yau's result is stronger in many respects,
most notably in the absence of any projectivity/rationality assumption and a more precise 
regularity theory. 
On the other hand, the condition on the poles of the L.H.S. is less optimal than here. 
We expect the projectivity/rationality assumptions to be superfluous. 

Observe that the condition on the singularity in the L.H.S. is precisely the condition
that the singular metric associated to $l(\tau_i)$ has a trivial multiplier ideal sheaf, or if 
$q=1$ that the pair $(X,l (\tau_1))$ is klt (see definition 6.7). 
The possibility of solving 
complex Monge-Amp\`ere equations with $L^p$-R.H.S. 
was first established by S.Kolodziej [K 1,2],
and the connection with
the singularities of the 
Minimal Model Program (MMP for short)
has been a strong incentive to our work. 
 
>From an algebraic geometer's perspective,
 these results may be viewed as a version of [Y] for normal K\"ahler spaces.

\vskip.2cm

As a by-product, S.T.Yau constructed K\"ahler-Einstein metrics on smooth
canonically polarized manifolds and Ricci-flat metrics on what is now known as 
Calabi-Yau manifolds. 
It had been soon realized [Ko] that this also yields K\"ahler-Einstein metrics 
on K\"ahler orbifolds hence on the
canonical models of surfaces of general type since they have isolated quotient singularities. 

In higher  dimension, in spite of the development of the MMP
during the 1980s -- culminating with [Mo] and the proof of the existence 
of canonical models for general type 3-folds [Ka] --, 
there was no satisfying analog of these K\"ahler-Einstein metrics. 

For smooth minimal general type projective manifolds, 
H.Tsuji proved [Ts] that an appropriate K\"ahler-Ricci flow starting 
with an arbitrary K\"ahler datum exists in infinite time,
converges towards a current representing the canonical
class which is smooth outside the exceptional locus of the map 
to the canonical model, and defines a K\"ahler-Einstein metric there 
\footnote{Although his idea was rather compelling,
the details of the proof for convergence were somewhat hard 
to follow.}. The conjecture made there that the current has 
continuous (or even bounded) local potentials partly
motivated our work.

The article [Ts] has been revisited in two recent preprints,
[CN] and [TZ], we learnt of when finishing the present work,
where a very satisfactory proof of convergence towards a current with bounded potentials is given. 
The independent work [TZ] uses a slightly weaker version of Theorem A 
and does not give any detail on the proof. These details have been indeed provided subsequently in [Zha]. 
On the other hand, the three approaches tend to emphasize different aspects of 
the problem and seem to complement each other nicely. 
  
In this article we give a more general theory of singular K\"ahler-Einstein metrics 
as a consequence of the following result:
\vskip.2cm
 
\noindent {\bf{Theorem C.}} 
{\em Let $(V,\Delta)$ be a projective klt pair such that $K_V+\Delta$ is an ample $\Q$-divisor.
  Then there is a unique semi-K\"ahler current in $[K_V+\Delta]$
  with bounded potentials,
  which satisfies a global degenerate Monge-Amp\`ere equation on $V$
  and defines a smooth
  K\"ahler-Einstein metric of negative curvature on $(V-\Delta)^{reg}$.
  
  Let $(V,\Delta)$ be a projective klt pair such that $K_V+\Delta\cong 0$ 
 ($\Q$-linear equivalence of $\Q$-Cartier divisors). Then in every ample class in $NS_{\R}(V)$
  there is a unique semi-K\"ahler current with bounded potentials,
  which satisfies a global degenerate Monge-Amp\`ere equation on $V$
  and defines a Ricci flat metric on $(V-\Delta)^{reg}$.} 
\vskip.2cm
  
The precise formulation and Monge-Amp\`ere equations are to be found in 
Theorems 7.5 and 7.8 below.

\vskip.2cm

 A large part of the MMP is now confirmed to work in higher dimension. Based on ideas of [Sho], developped further on by [HMcK], 
 the preprint [BCHM] has achieved a full proof of
 the finite generation of the canonical ring for any projective manifold.   Using this, we obtain:
 
 \vskip.2cm
  
\noindent {\bf{Corollary D.}}
{\em Let $X$ be a projective manifold of general type and $V=X_{can}$ the unique
  model of $X$ with only canonical singularities and $K_V$ ample. Then $K_V$ contains a unique 
  singular K\"ahler-Einstein metric $\omega_{KE}$ of negative curvature. }
  
\vskip.2cm
  
Note that we do not assume the singularities to be quotient singularities nor 
that $X$ has a smooth minimal model
(a strong restriction present in [Ts], [TZ], [CN]). On the other hand if $\pi: X\to X_{can}$ 
is a resolution of singularities then $[\pi^* \omega_{KE}] \cong K_X +F$ where $F\ge 0$ 
and $=0$ iff $\pi$ is crepant and $X$ is a smooth minimal model.

  The problem of constructing a singular K\"ahler-Einstein metric on a canonically polarized projective variety 
with canonical singularities $X_{can}$ had already been considered in [Sg]\footnote{We apologize for having overlooked
 this reference in the  first circulated version of this work. }, with a different approach. Theorem 5.6 there and its proof imply that given
$\pi: X\to X_{can}$ a log resolution, there is a closed positive current $T_{KE}$ in $\pi^*K_{X_{can}}$ 
with zero Lelong numbers such that $T_{KE}$ is smooth on $\pi^{-1} (X_{can}^{reg})$ and defines a KE metric there.
This constuction agrees with ours, and our contribution is that  $T_{KE}$ has a bounded potential coming from $X_{can}$.

One may also try and construct these singular KE metrics unconditionnally
in an approach to the finiteness of the canonical ring [Siu 2]. 
Needless to say, it is a substantially harder task.

\vskip.2cm

Let us state yet another corollary:

\vskip.2cm

\noindent {\bf{Corollary E.}}
{\em Let $X$ be a projective variety with only canonical singularities such that $K_X \sim_{{\mathbb Q}} 0$ and $A\in NS_{{\mathbb R}}(X)$ an ample class. Then $A$ contains a unique singular Ricci-flat  metric $\omega_{CY}$ with bounded potentials
  
  The local potentials of $\omega_{CY}$ are continuous provided $(V,A)$ satisfies the continuous approximation property. }
\vskip.2cm

For the applications, the regularity theory of singular K\"ahler-Einstein metrics still needs to be developped more thoroughly. 

\section*{Notations and organization of the paper}

In the whole paper, $X$ will denote a compact K\"ahler manifold
of dimension $n$, 
$\o$ a smooth closed form of bidegree $(1,1)$ which
is non-negative and {\it big}, i.e. the smooth measure $\o^n$ is not identical to zero.
For convenience
we normalize $\o$ so that
$$
Vol_{\o}(X):=\int_X \o^n =1.
$$
 
 $V$ will denote a
normal complex space. 
A resolution of $V$ will be a locally projective bimeromorphic holomorphic morphism $\pi:X\to V$, 
$X$ being smooth, such that 
$\pi: \pi^{-1} (V^{reg})\to V^{reg}$ is an isomorphism. A resolution $\pi$ is a log-resolution
iff $\pi^{-1} (V^{sing})$ is a divisor with simple normal crossings.
Assume we have a coherent ideal sheaf $\mathcal I \subset \Ox$. A log resolution of
$(V,\mathcal {I})$ is a projective bimeromorphic holomorphic morphism $\pi: X\to V$
$X$ being smooth, such that 
$\pi: \pi^{-1} (V-Z(\mathcal{I}))^{reg})\to (V-Z(\mathcal{I}))^{reg}$ is an isomorphism 
with the additional property that the ideal sheaf $\pi^{-1}\mathcal{I} . \mathcal{O}_X$
 \footnote{ $\pi^{-1}\mathcal{I} . \mathcal{O}_X$ is locallly the ideal sheaf of 
$\mathcal{O}_X$ generated by the family of holomorphic functions
 $(\pi^* f_I)_I$, where $(f_I)_I$ are local generators of $\mathcal {I}$. 
The set $Z(\mathcal{I})$ is the analytic subvariety defined by $\mathcal{I}$. } 
satisfies  $\pi^{-1}\mathcal{I} . \mathcal{O}_X = \mathcal{O}_X (-\sum \gamma_E E)\subset \mathcal{O}_X$
 where $\gamma_E\in \N$ is a positive integer attached to an exceptional divisor $E$ of $\pi$. 
 
 A pair is 
a pair $(V,\Delta)$ with $V$ a normal complex space and $\Delta$ a $\Q$-Weil divisor
$\Delta=\sum_i d_iE_i$ where $0\le d_i\le 1$ are rational numbers and $(E_i)_i$ is a finite family
of pairwise distinct irreducible codimension 1 subvarieties of $V$.
A log resolution of a pair is a log resolution of the ideal $\mathcal{I}_{N\Delta}$
 where $N$ is an integer such that $Nd_i\in \N$. 
\footnote{The MMP is conjectured to work for $\mathbb{Q}$-factorial dlt pairs provided $V$ is projective algebraic. This seemingly technical 
extension of the MMP is known as log-MMP. log-MMP works in dimension $\le 3$. 
The philosophy of the log-MMP is to define the canonical divisor of a pair to be
 $K_{(V,\Delta)}:=K_V+\Delta$ and to try and prove the same theorems for pairs and 
for varieties. See [BCHM] for the strongest results in this direction.}
 
All these flavors of log-resolutions exist by [Hi], [BM] if the
 variety (resp. pair) under consideration is
open in (resp. a restriction to an open subset of) a compact variety (pair).
The resolution can then be assumed to be a projective morphism.

\vskip.2cm

The paper is organized as follows.
In {\it section 1} we define, following [GZ 2],
the set ${\mathcal E}^1(X,\o)$ of $\o$-psh functions with finite self-energy,
and produce weak solutions to complex Monge-Amp\`ere equations
$(\o+dd^c \f)^n=\mu$ in the class ${\mathcal E}^1(X,\o)$
(see proposition 1.4). This is our first basic observation: weak solutions are easy
to produce in ${\mathcal E}^1(X,\o)$.

The continuity of the solutions is studied in {\it section 2}
(see Theorem 2.1), by using ideas from [K 1,2,3] and [GZ 1].
This, together with propositions 3.1 and 3.3, yields Theorem A.
We actually expect the solutions to be H\"older-continuous, 
as Theorem 3.5 indicates.
We indeed establish further regularity results in
{\it section 3}, especially Theorem 3.6, by using ideas of [Y],[Ts].
This yields Theorem B.

In {\it section 4} we solve Monge-Amp\`ere equations of the
type $(\o+dd^c \f)^n=e^{t \f} \mu$, $t>0$
(see Theorems 4.1, 4.4) by a fixed point method.
Here again the use of class ${\mathcal E}^1(X,\o)$ makes life
easier, and allows us to reduce our analysis to previously studied 
Monge-Amp\`ere equations $(\o+dd^c \f)^n=\mu'$.

In {\it section 5} we recall some basic facts on some of the singularities
encountered in the MMP, and in {\it section 6} we explain what sort of
measures $\mu$ we need to consider in order to produce K\"ahler-Einstein metrics.
An important observation is lemma 6.4, which shows that it is necessary
to restrict to the case of log terminal singularities (see definition 5.3).

In {\it section 7} we show how our results from sections 2,3,4 allow
us to produce singular K\"ahler-Einstein metrics
(see Theorems 7.5, 7.8, 7.12). This is where we prove Theorem C.

\section{Weak solutions to Monge-Amp\`ere equations}

In this section $\o$ will denote
a smooth semi-positive closed $(1,1)$-form which is big, i.e. satisfies $\int_X\omega^n>0$. 
For any K\"ahler form $\Omega$ on $X$ and for all $\e>0$,
the form $\o_{\e}:=\o+\e \Omega$ is again K\"ahler on $X$.

We are going to extend several results that are known to hold true 
when $\o$ is K\"ahler to the present more general setting.

Recall that the set of $\o$-plurisubharmonic functions ($\o$-psh for short)
is
$$
PSH(X,\o):=\{ \f \in L^1(X,\R \cup \{ -\infty \}) \, / \, 
dd^c \f \geq -\o \text{ and } \f
\text{ is } \text{u.s.c.} \}.
$$
We refer the reader to [GZ 1] for basic properties of $\o$-psh functions.
The following subclass has been extensively studied in [GZ 2]:

\begin{defi}
We let ${\mathcal E}^1(X,\o)$ denote the set of $\o$-psh functions with finite self-energy:
this is the set of functions $\f \in PSH(X,\o)$ for which
there exists a sequence $\f_j \in PSH(X,\o) \cap L^{\infty}(X)$
such that
$$
\f_j \searrow \f \; \; 
\text{ and } \; \; 
\sup_j \int_X (-\f_j) (\o +dd^c \f_j)^n <+\infty .
$$
\end{defi}

This class of functions is studied in [GZ 2] when $\o$ is a K\"ahler form.
We leave it to the reader to check that the basic properties
of this class of functions proved in [GZ 2] when $\o$ is K\"ahler
apply with no modification to the present case. 
In particular the complex Monge-Amp\`ere operator
$(\o+dd^c \f)^n$ is well-defined for $\f \in {\mathcal E}^1(X,\o)$,
and it is continuous on decreasing sequences of functions in ${\mathcal E}^1(X,\o)$.

We shall need a slightly more general continuity result, which takes into
account the dependence in $\o$:

\begin{pro}
Fix $\Omega$ a K\"ahler form on $X$, and
let $(\e_j)$ be a sequence of positive real numbers decreasing to zero.
Let $\f_j \in {\mathcal E}^1(X,\o +\e_j \Omega)$ be a sequence of functions which 
decrease pointwise towards $\f$, and such that
$$
\sup_{j \geq 1} \int_X |\f_j| (\o+\e_j \Omega+dd^c \f_j)^n <+\infty.
$$

Then $\f \in {\mathcal E}^1(X,\o)$, and 
$(\o+\e_j \Omega+dd^c \f_j)^n \rightarrow (\o+dd^c \f)^n$.
\end{pro}

\begin{proof}
Set $\o_j:=\o+\e_j \Omega$. We can assume w.l.o.g.
that $\f \leq \f_j \leq 0$. Set
$$
\f^K:=\max ( \f, -K) \in PSH(X,\o)
\; \text{ and } \; 
\f_j^K:= \max (\f_j,-K) \in PSH(X,\o_j).
$$
Observe that, $K$ being fixed, $(\f_j^K)_j$ is uniformly bounded and decreases
towards $\f^K$ as $j$ goes to infinity. Therefore
$(\o_j+dd^c \f_j^K)^n \rightarrow (\o+dd^c \f^K)^n$, by a classical result
of E.Bedford and A.Taylor [BT 82].
Moreover the sequence of positive measures $(-\f_j^K)(\o_j+dd^c \f_j^K)^n$
has uniformly bounded mass, since by lemma 2.3 in [GZ 2],
$$
0 \leq \int_X (-\f_j^K)(\o_j+dd^c \f_j^K)^n
\leq 2^n \int_X  (-\f_j)(\o_j+dd^c \f_j)^n \leq 2^n M,
$$
where $M:=\sup_j \int_X  (-\f_j)(\o_j+dd^c \f_j)^n<+\infty$. 

Since $\f_j$ is u.s.c., a standard argument yields that any cluster point
$\nu$ of the sequence $(-\f_j^K)(\o_j+dd^c \f_j^K)^n$ satisfies
$0 \leq (-\f^K)(\o+dd^c \f^K)^n \leq \nu$. In particular
$$
0 \leq \int_X (-\f^K)(\o+dd^c \f^K)^n \leq 
\liminf_{j \rightarrow +\infty} \int (-\f_j^K) (\o_j+dd^c \f_j^K)^n
\leq 2^n M
$$
is bounded from above uniformly with respect to $K$. Since
$\f^K$ decreases towards $\f$, this shows $\f \in {\mathcal E}^1(X,\o)$.

It remains to show that $(\o_j+dd^c \f_j)^n \rightarrow (\o+dd^c \f)^n$.
Since $(\o_j+dd^c \f_j^K)^n \rightarrow (\o+dd^c \f^K)^n$ for any fixed $K$,
it is enough to get an upper bound on the mass of
$(\o_j+dd^c \f_j^K)^n$ in $(\f_j \leq -K)$ which is uniform in $j$.
This follows from Chebyshev inequality, namely
$$
\int_{(\f_j \leq -K)} (\o_j+dd^c \f_j)^n  \leq \frac{1}{K} 
\int_X (-\f_j) (\o_j+dd^c \f_j)^n \leq \frac{M}{K}.
$$
This yields the desired result.
\end{proof}
\vskip.2cm

The Monge-Amp\`ere capacity $Cap_{\o}(\cdot)$ has been studied in [GZ 1],
$$
Cap_{\o}(K):=\sup \left\{ \int_K \o_{u}^n \, / \, u \in PSH(X,\o), 
0 \leq u \leq 1 \right\},
$$
where $K$ is a Borel subset of $X$.
Here -- and in the sequel -- we use the notation $\o_u:=\o+dd^c u \geq 0$.
In this article we are interested in measures which are
dominated by the Monge-Amp\`ere capacity in the following way:

\begin{defi}
A probability measure $\mu$ on $X$ satisfies condition
${\mathcal H}(\a,A,\o)$ if for all Borel subset $K$ of $X$,
$$
\mu(K) \leq A Cap_{\o}(K)^{1+\a}.
$$
\end{defi}

It has been shown by S.Kolodziej that 
when $\o$ is {\it K\"ahler},
a probability measure $\mu$ which 
satisfies ${\mathcal H}(\a,A,\o)$ can be written as the Monge-Amp\`ere
measure of some continuous $\o$-psh function.
This is still true when $\o$ is merely {\it semi-positive} and {\it big},
and the proof will occupy us until the end of section 2.
We start by observing  -- following [GZ 2] -- that $\mu$ is the Monge-Amp\`ere 
of a function $\f$ which is not too singular.

\begin{pro}
Let $\mu$ be a probability measure on $X$ which satisfies condition
${\mathcal H}(\a,A,\o)$. 
Then there exists a unique function
$\f \in {\mathcal E}^1(X,\o)$ s.t.
$$
\mu=(\o+dd^c \f)^n \; \text{ and } \;
\sup_X \f=-1.
$$
\end{pro}

\begin{proof}
Fix $\Omega$ a K\"ahler form on $X$, and set 
$\o_j:=\o+\e_j \Omega$, where $\e_j>0$ decreases to $0$.
We start by showing that ${\mathcal E}^1(X,\o_j) \subset L^1(\mu)$.

Fix $\f \in {\mathcal E}^1(X,\o_j)$. We can assume without loss of
generality that $\sup_X \f = -1$. 
It follows from propositions 3.6 and 2.7 in [GZ 1]
that there exists a constant $C=C(\o,\Omega)>0$ independent of $j$ such that
$Cap_{\o_j}(\f<-t) \leq C/t$ for all $t>0$. 
Since $Cap_{\o}(\cdot) \leq Cap_{\o_j}(\cdot)$, the measure $\mu$ satisfies
${\mathcal H}(\a,A,\o_j)$.
We infer
\begin{equation}
0 \leq \int_X (-\f) d\mu =\int_{t=1}^{+\infty} \mu(\f<-t)dt
\leq \frac{A C^{1+\a}}{\a} <+\infty,
\end{equation}
with an upper-bound which is independent of $j$.

The main result in [GZ 2] guarantees in this case that there exists a unique
function $\f_j \in {\mathcal E}^1(X,\o_j)$ such that
$$
(\o_j+dd^c \f_j)^n=\l_j \mu \, \text{ and } \; \sup_X \f_j=-1,
$$
where $\l_j=\int_X (\o+\e_j \Omega)^n >1$ decreases to 1 as $j$ goes to infinity.

The normalization $\sup_X \f_j=-1$ implies that the sequence $(\f_j)$ is relatively
compact in $L^1(X)$ (see proposition 2.7 in [GZ 1]).
Let $\f$ be a cluster point of $(\f_j)$. Relabelling if necessary, we 
assume $\f_j \rightarrow \f$ in $L^1(X)$. Note that $\f \in PSH(X,\o)$
and $\sup_X \f=-1$ (by Hartogs' lemma, see proposition 2.7, [GZ 1]).
We are going to show that $\f \in {\mathcal E}^1(X,\o)$ and $\o_{\f}^n=\mu$.

Set $\Phi_j:=(\sup_{l \geq j} \f_l)^*$, where $u^*$ denotes the upper-semi-continuous 
regularization of $u$. Then $\Phi_j \in PSH(X,\o_j)$ with $\Phi_j \geq \f_j$,
hence $\Phi_j \in {\mathcal E}^1(X,\o_j)$
(see proposition 3.2 in [GZ 2]), and $\Phi_j$ decreases towards $\f$.
For $l \geq j$, we have
$$
(\o_j+dd^c \f_l)^n \geq (\o_l+dd^c \f_l)^n=\l_l \mu \geq \mu.
$$
It follows therefore from an inequality due to J.-P.Demailly [Dem 1]  that
$(\o_j+dd^c \Phi_j)^n \geq \mu$. Now by (1) and lemma 2.3 in [GZ 2],
$$
0 \leq \int_X (-\Phi_j) (\o_j+dd^c \Phi_j)^n 
\leq 2^n \int_X (-\f_j) (\o_j+dd^c \f_j)^n
=2^n \l_j \int_X (-\f_j) d\mu ,
$$
is uniformly bounded with respect to $j$ thanks to (1).

We infer from proposition 1.2 that $\f \in {\mathcal E}^1(X,\o)$
and $(\o_j+dd^c \Phi_j)^n \rightarrow (\o+dd^c \f)^n$.
Thus $(\o+dd^c \f)^n \geq \mu$, but these are two probability measures,
whence $\mu=(\o+dd^c \f)^n$. The uniqueness of $\f$ follows from Theorem 3.3, [GZ 2].
\end{proof}

\section{Continuous solutions}

The goal of this section is to prove the following result:

\begin{thm}\label{mas}
Let $\mu$ be a probability measure on $X$ which satisfies condition
${\mathcal H}(\a,A,\o)$. Then there exists a unique
{\it bounded} function $\f \in PSH(X,\o)$ such that
$$
\mu=(\o+dd^c \f)^n \; \text{ and } \; \sup_X \f=-1.
$$
Moreover $||\f||_{L^{\infty}(X)} \leq C$, where $C$ only
depends on $\a,A$ and $\o$.  

Furthermore $\f$ is continuous provided $(X,[\omega])$ satisfies the continuous approximation property.
\end{thm}

\begin{defi} Let $X$ be  a compact K\"ahler manifold and $\omega$ a smooth semipositive closed $(1,1)$-form. We say $(X,[\omega])$ satisfies the continuous approximation property  if for every $\omega$-psh bounded function $\f$ there is a decreasing sequence of
continuous $\omega$-psh functions such that $\lim \f_j=\f$ pointwise. 
\end{defi}

It follows from Demailly's regularisation theorem that if $X$ is smooth and $[\omega]$ K\"ahler then $(X,[\omega])$ satisfies the continuous approximation property, see [BK] for a simple proof. 

The authors believe the approximation property is satisfied if $X$ is a resolution of singularities of a normal projective variety with canonical singularities and $[\omega]\in NS (X)$ is the pull-back of an ample $A$
 class of $V$. If so, one says that $(V,A)$ enjoys the continuous approximation property.

It follows from proposition 1.4 that we already know the
existence of a unique solution $\f \in {\mathcal E}^1(X,\o)$
to this Monge-Amp\`ere equation. We need to show it is continuous.
The following result is the key to everything to follow.

\begin{lem}
Let $\f,\p \in {\mathcal E}^1(X,\o)$ be two negative functions.
Then for all $s>0$ and $0 \leq t \leq 1$,
$$
t^n Cap_{\o}(\f-\p<-s-t) \leq \int_{(\f-\p<-s-t\p)} \o_{\f}^n.
$$
\end{lem}

\begin{proof}
Fix $u \in PSH(X,\o)$ with $0 \leq u \leq 1$. 
For $\d>0$ we set $t=\d/(1+\d)$.
Observe that $0 \leq t \leq 1$ and
$$
\{ \f-\p<-s-t\} \subset \left\{\f < \frac{\p+\d u}{1+\d} -s-t \right\} 
\subset \{ \f -\p<-s-t\p\}.
$$
Set $\tilde{\f}:=(\p+\d u)/(1+\d) -s-t \in PSH(X,\o)$.
Observe that
$$
t^n \int_{(\f-\p<-s-t)} \o_u^n 
\leq \int_{(\f-\p<-s-t)} \left[ \frac{1}{1+\d} \o_{\p}+\frac{\d}{1+\d} \o_u \right]^n
\leq \int_{(\f<\tilde{\f})} [\o+dd^c \tilde{\f}]^n.
$$
It follows from the
comparison principle in class ${\mathcal E}^1(X,\o)$,
that
$$
\int_{(\f<\tilde{\f})} [\o+dd^c \tilde{\f}]^n
\leq \int_{(\f <\tilde{\f})} [\o+dd^c \f]^n
\leq \int_{(\f-\p<-s-t\p)} [\o+dd^c \f]^n.
$$
Taking the supremum over all u's yields the desired result.
\end{proof}

We will also need the following elementary observation:

\begin{lem}
Let $f:\R^+ \rightarrow \R^+$ be a decreasing right-continuous function 
such that $\lim_{+\infty } f=0$. Assume 
there exists $\a,B>0$ such that
$f$ satisfies
$$
H(\a,B) \hskip1cm
tf(s+t) \leq B [ f(s) ]^{1+\a}, \;
\forall s>0, \, \forall 0 \leq t \leq 1.
$$

Then there exists $S_{\infty}=S_{\infty}(\a,B) \in \R^+$ such that
$f(s)=0$ for all $s \geq S_{\infty}$.
\end{lem}

\begin{proof}
Fix $s_0 >0$ large enough so that $f(s_0)^{\a}<1/2B$. We define a sequence
$(s_j) \in \R_+^{\N}$ by induction in the following way.
If $f(s_0)=0$ we stop here, otherwise we set
$$
s_1:=\sup \left\{ s>s_0 \, / \, f(s)>\frac{1}{2} f(s_0) \right\}.
$$
Observe that $s_1 \leq 1+s_0$ thanks to $H(\a,B)$ and by definition of $s_0$.

Since $f$ is right-continuous we get $f(s_1) \leq f(s_0)/2$. If $f(s_1)=0$
we stop here, otherwise we go on by induction, setting
$$
s_{j+1}:=\sup \left\{ s>s_j \, / \, f(s)>\frac{1}{2} f(s_j) \right\}.
$$
At each step $f(s_{j+1}) \leq f(s_j)/2$ and $s_{j+1} \leq 1+s_j$.
However the sequence $(s_j)$ does not grow too fast. It follows indeed
from $H(\a,B)$ that if $s \in ]s_j,s_{j+1}[$,
$$
(s-s_j) f(s) \leq B f(s_j)^{1+\a} \leq 2B f(s) f(s_j)^{\a},
$$
since $f(s_j)/2 \leq f(s) \leq f(s_j)$ on the interval $[s_j,s_{j+1}]$.
We infer
$$
s_{j+1}-s_j \leq 2Bf(s_j)^{\a} \leq 2B 2^{-j\a}f(s_0)^{\a} \leq 2^{-j\a}.
$$
Thus the sequence $(s_j)$ is bounded from above, with limit
$$
S_{\infty}=s_0+\sum_{j \geq 0}(s_{j+1}-s_j) \leq
s_0+\frac{2B f(s_0)^{\a}}{1-2^{-\a}} \leq s_0+\frac{1}{1-2^{-\a}}.
$$
\end{proof}

\begin{rqes}
Observe that the starting time $s_0(f,\a,B)$ is invariant under 
dilatation $f \mapsto \l f$, which transforms $B$ into $B/\l^{\a}$.
Note also that if $f(0)^{\a}<1/2B$, then we can take $s_0=0$, hence we get
in this case 
$$
S_{\infty} \leq \frac{2B}{1-2^{-\a}} [f(0)]^{\a}.
$$
\end{rqes}

To see how previous lemmas can be used, we first prove that
the unique solution $\f \in {\mathcal E}^1(X,\o)$ given by
proposition 1.4 is bounded.
\vskip.3cm

\noindent {\it A uniform bound on the solution.}
Let $\f \in {\mathcal E}^1(X,\o)$ be the unique function
such that $\mu=(\o+dd^c \f)^n$ and $\sup_X \f=-1$.
Set
$$
f(s):=[Cap_{\o}(\f<-s)]^{1/n}.
$$
Observe that $f:\R^+ \rightarrow \R^+$ is a right-continuous decreasing
function with $\lim_{+\infty} f=0$.
Since $\mu=\o_{\f}^n$ satisfies ${\mathcal H}(\a,A,\o)$, it follows from
lemma 2.2 applied to the function $\p \equiv 0$ that
$f$ satisfies $H(\a,B)$ with $B=A^{1/n}$.

It follows from propositions 2.7 and 3.6 in [GZ 1] that
$f(s) \leq C_1/s^{1/n}$ for some constant $C_1$ which only depends on $\o$.
We can thus take a starting time $s_0=2^{n/\a}C_1^nA^{1/\a}$ (see the
proof of lemma 2.3) and get $f(s)=0$ for
$s \geq S_{\infty}:=s_0+(1-2^{-\a})^{-1}$: this shows that the sets
$(\f <-s)$ are empty if $s \geq S_{\infty}$, hence
\begin{equation}
||\f||_{L^{\infty}(X)} \leq 2^{n/\a} C_1^n A^{1/\a}+\frac{1}{1-2^{-\a}}.
\end{equation}
\vskip.3cm

We are now going to use a refinement of the previous reasoning 
in order to show that $\f$ is actually continuous if the continuous approximation property holds.

\begin{pro}
Let $\f,\p \in {\mathcal E}^1(X,\o)$ be two  functions such that 
$\sup_X \f= \sup_X \p =-1$ 
and fix  $\e>0$. 
Assume $\o_{\f}^n=\mu$ satisfies ${\mathcal H}(\a,A,\o)$
and $\p$ is bounded.
There exists $C=C(\a,A,\o, ||\p||_{L^{\infty}(X)})>0$ such that
$$
\sup_X (\p-\f) \leq \e+C \left[ Cap_{\o}(\f-\p<-\e) \right]^{\a/n}.
$$
\end{pro}

This inequality can be interpreted as follows.
Assume $\p$ is also a solution to a Monge-Amp\`ere equation
$\o_{\p}^n=\mu'$, where $\mu'$ also satisfies ${\mathcal H}(\a,A,\o)$.
Then we can interchange the roles of $\f$ and $\p$ and
get an upper-bound on $||\f-\p||_{L^{\infty}(X)}$. 
The proposition then tells us that if $\f$ and $\p$ are close
in capacity, they are close in $L^{\infty}$-norm.

\begin{proof}
Set $M:=||\p||_{L^{\infty}(X)}$. Observe that when $t \geq 0$,
$(\f-\p<-s-t\p) \subset (\f-\p<-s+tM)$, hence it follows from lemma 2.2
that for all $s>0$, $0 \leq t \leq 1$,
\begin{equation}
t^n Cap_{\o}(\f-\p<-s-t) \leq (1+M)^n \int_{(\f-\p<-s)} \o_{\f}^n,
\end{equation}
using the obvious substitutions $s \mapsto s-Mt$, $t \mapsto (1+M)t$.
Since $\mu=\o_{\f}^n$ satisfies ${\mathcal H}(\a,A,\o)$, we infer
$$
t^n Cap_{\o}(\f-\p<-s-t) \leq A (1+M)^n Cap_{\o}(\f-\p<-s)^{1+\a}.
$$
Consider
$$
f(s):=\left[ Cap_{\o}(\f-\p<-s-\e)\right]^{1/n}, \, s>0.
$$
Then $f$ satisfies the condition $H(\a,B)$ of lemma 2.3 with
$B=(1+M) A^{1/n}$.
Assume
$f(0)=\left[ Cap_{\o}(\f-\p<-\e)\right]^{1/n}<\frac{1}{(2B)^{1/\a}}$.
It follows in this case from Remarks 2.4
that $f(s)=0$ for $s \geq S_{\infty}$, where
$$
S_{\infty} \leq \frac{2B}{1-2^{-\a}} \left[ Cap_{\o}(\f-\p<-\e)\right]^{\a/n}
$$
Therefore the sets $\{ \f-\p <-s-\e \}$ are empty for
$s>S_{\infty}$, hence
$$
\sup_X (\p-\f) \leq \e+S_{\infty} \leq \e+
\frac{2(1+M)A^{1/n}}{1-2^{-\a}} \left[ Cap_{\o}(\f-\p<-\e)\right]^{\a/n}.
$$
Thus we can take here $C \geq 2B/(1-2^{-\a})$.

If $f(0)=[Cap_{\o}(\f-\p<-\e)]^{1/n} \geq (2B)^{-1/\a}$, then
$$
\sup_X (\p-\f) \leq -\inf_X \f \leq C_2(\a,A,\o),
$$
by (2), hence it suffices to take $C \geq 2B C_2 $ to conclude.
\end{proof}

\begin{proof2.1}
Let $\f \in {\mathcal E}^1(X,\o)$ be the unique solution to
the normalized Monge-Amp\`ere equation
$\mu=(\o+dd^c \f)^n$, $\sup_X \f=-1$ (see proposition 1.4).
It follows from (2) that $\f$ is bounded.

If $(X,[\omega])$ enjoys the continuous approximation property, let $\f_j$ be 
decreasing sequence of $\omega$-psh functions decreasing pointwise to  $\f$.
Since $\sup_X \f=-1$, we can assume $\f_j \leq 0$.
Since $\f$ is bounded and $\f_j \geq \f$, the functions $\f_j$ are uniformly
bounded on $X$.

Observe that $\f_j$ converges towards $\f$ in capacity
(see proposition 3.7 in [GZ 1]), hence $\lim Cap_{\o}(\f-\f_j<-\e) =0$,
for all $\e>0$.
It follows therefore from proposition 2.5 (applied with $\p=\f_j$) that for all $\e>0$,
$$
\lim_{j  \rightarrow +\infty} ||\f-\f_j||_{L^{\infty}(X)}
=\lim_{j  \rightarrow +\infty} \sup_X (\f_j-\f)
\leq \e.
$$
Thus $(\f_j)$ converges uniformly towards $\f$, hence 
$\f$ is continuous.
\hfill
\end{proof2.1}

\section{More regularity}

\subsection{Measures with density}

We now turn to the study of the complex Monge-Amp\`ere equation
$$
(\o+dd^c \f)^n=\mu, 
\text{ when } \mu=f \o^n
$$
is a measure with
density $0 \leq f \in L^p(\o^n)$, $p>1$.

\begin{pro}\label{lp}
Assume $\mu=f \o^n$ is a probability measure with density
$0 \leq f \in L^p(X)$, for some $p>1$.
Then for any $\a>0$, there exists $A_{\a}>0$ such that
$\mu$ satisfies ${\mathcal H}(\a,A_{\a},\o)$.
\end{pro}

\begin{proof}
It is enough to establish ${\mathcal H}(\a,A_{\a},\o)$ for compact subsets,
by regularity of $\mu$ and $Cap_{\o}$.
Let $K$ be a compact subset of $X$. It follows from H\"older's inequality
that
$$
0 \leq \mu(K) \leq ||f||_{L^p(\o^n)} \left[ \text{Vol}_{\o}(K)\right]^{1/q},
$$
where $1/p+1/q=1$. 
Note that $||f||_{L^p(\o^n)}=1$ since we assume $\mu$ is a probability measure.
We claim that
\begin{equation}
\text{Vol}_{\o}(K) \leq C_{\o} \exp \left[ -\g_{\o}(Cap_{\o}(K))^{-1/n} \right],
\end{equation}
for some constants $C_{\o},\g_{\o}>0$ that only depend on $\o$.
We will be done if we can prove (4) since we can then check by elementary computations
that $\exp(-x^{-\d})$ is dominated from above by
$A_{\a} x^{\a}$, for all $x \in [0,1]$.

The set of functions ${\mathcal F}_0:=\{ \f \in PSH(X,\o) \, / \, \sup_X \f=0 \}$
is compact in $L^1(X)$ (see proposition 2.7, [GZ 1]). These functions have Lelong
numbers $\nu(\f,x) \leq \nu_{\o}$ bounded from above by a uniform constant.
It follows therefore from Skoda's uniform integrability theorem [Z],
that
$$
\sup_{\f \in {\mathcal F}_0} \int \exp \left[ -\frac{2 \f}{\nu_{\o}+1} \right] \o^n 
\leq C_2 <+\infty.
$$
Set $\g_{\o}:=2/(\nu_{\o}+1)>0$ and let
$$
V_{K,\o}^*(x):= \left( \sup \{ \f(x) \, / \, \f \in PSH(X,\o), \, 
\f \leq 0 \text{ on } K \} \right)^*
$$
denote the Siciak extremal function of $K$
(see section 5.1 in [GZ 1]).
Then
$$
Vol_{\o}(K) \leq \int_X \exp \left( -\g_{\o} V_{K,\o}^* \right) \o^n
\leq C_2 T_{\o}(K)^{\g_{\o}},
$$
where $T_{\o}(K):=\exp(-\sup_X V_{K,\o}^*)$ denotes the Alexander capacity
of $K$.
It follows now from theorem 7.1 in [GZ 1] that
$$
T_{\o}(K) \leq e \exp \left[ -Cap_{\o}(K)^{-1/n} \right],
$$
which yields (4).
\end{proof}

It follows therefore from theorem 2.1 that there exists a unique continuous function
$\f \in PSH(X,\o)$ such that
$$
\mu=f\o^n=(\o+dd^c \f)^n, \text{ with } \sup_X \f=-1,
$$
when $0 \leq f \in L^p(\o^n)$, $p>1$, with $\int_X f \o^n =1$.

Actually we will be interested in measures
with $L^p$-density with respect to a positive definite volume form $d\l$,
while the smooth measure $\o^n$ may vanish along a divisor.
This does not make much difference, as follows from H\"older's inequality:

\begin{lem} \label{leb} Let $V$ be a $n$-dimensional compact normal K\"ahler space 
and $\Omega$ be a smooth 
K\"ahler form on $V$. Let $\pi: X\to V$ a resolution, $\o=\pi^* \Omega$, 
and let $d\lambda$ be a positive definite smooth
volume form on $X$. 

If $\mu=f_1d\lambda$, with $f_1\in L^p(X, d\lambda)$ for some $p>1$,
then there exists $p'>1$ such that $\mu= f \omega^n$ and $f \in L^{p'} (X,\omega^n)$. 
\end{lem}

\begin{proof}
Observe that $\omega^n = E d\lambda$ for some smooth density $E \geq 0$ 
which vanishes along the exceptional locus of $\pi$,
thus
$$
\mu= f_1 d\l =f \o^n,
\text{ where } f =f_1/E.
$$

Fix local coordinates $(z^i)_{1\le i\le n}$ on a polydisk $\mathbb{D}\subset X$ and a local 
embedding $F:V\to \C^m$.
Note that $E$ is comparable to 
$\left|\frac{\partial F}{\partial z^1} \wedge \ldots \wedge  \frac{\partial F}{\partial z^1} \right|^2
\simeq \sum_{i=1}^r |f_i|^2$, $f_i$ being holomorphic on $\mathbb{D}$. 
Therefore $E\in L^{\infty}_{loc} (\mathbb{D})$
and $E^{-\a} \in L_{loc}^1 (\mathbb {D}, d\lambda)$ for
some $0 <\a<1$.

Choose $0<\alpha'<\alpha$ such that $\frac{1}{p} + \frac{1}{\alpha} =\frac{1}{\alpha'}$. 
Then $f^{\alpha'}=
f_1^{\alpha'} E^{-\alpha'}$ is the product of a function in $L^{p/\alpha'}(d\lambda)$
and a function in $L^{\alpha/ \alpha'}(d\l)$, hence it is in $L^1_{loc} (\mathbb{D}, d\lambda)$ 
by H\"older's inequality. 
A second application of H\"older's inequality yields
$$
\int_{{\mathbb{D}}} f^{1+\e} \o^n=
\int_{\mathbb{D}} f^{\epsilon} f_1 d\lambda 
\leq
\left(\int_{\mathbb{D}} f^{\epsilon q} d\lambda \right)^{1/q} 
\left(\int_{\mathbb{D}} f_1^p d\lambda \right)^{1/p}
<+\infty,
$$
where $q$ denotes the conjugate exponent to $p$.
This shows that $f \in L^{p'}(\o^n)$ if $p'=1+\e>1$ is chosen so small
that $\e q<\a'$.
\end{proof}

\subsection{H\"older continuity}

\begin{pro}
Assume $\o_{\f}^n=f \o^n$, $\o_{\p}=g \o^n$, where 
$\f,\p \in PSH(X,\o)$ are bounded and 
$f,g \in L^p(\o^n)$, $p>1$.
Then for all $0 < \g < 2/(2+nq)$,
$$
||\f-\p||_{L^{\infty}(X)} \leq C ||\f-\p||_{L^2(\o^n)}^{\g},
$$
where $q=p/(p-1)$ denotes the conjugate exponent to $p$.
\end{pro}

\begin{proof}
Fix $\e>0$ and $\a>0$ to be chosen later.
It follows from (2) and propositions 2.5, 3.1 that
$$
||\f-\p||_{L^{\infty}(X)} \leq \e+ C_1 \left[ Cap_{\o}(|\f-\p|>\e) \right]^{\a/n}.
$$
Applying the refined version of lemma 2.2 which involves
the uniform bound on $||\f||_{L^{\infty}(X)}, ||\p||_{L^{\infty}(X)}$
(see inequality (3)), we obtain
$$
Cap_{\o}(|\f-\p| >\e) \leq \frac{C_2}{\e^{n+2/q}} \int_X |\f-\p|^{2/q} (f+g) \o^n.
$$
It follows thus from H\"older's inequality that
$$
Cap_{\o}(|\f-\p|>\e) \leq \frac{C_3 ||f+g||_{L^p}}{\e^{n+2/q}}
\left[ ||\f-\p||_{L^2(\o^n)} \right]^{2/q}.
$$
Choose now $\e:=||\f-\p||_{L^2}^{\gamma}$ where
$0 < \g <2/(2+nq)$. Then
$$
Cap_{\o}(|\f-\p|>\e) \leq C_4 \left[ ||\f-\p||_{L^2} 
\right]^{2/q-\g(n+2/q)}.
$$
We infer
$$
||\f-\p||_{L^{\infty}(X)} \leq ||\f-\p||_{L^2}^{\g}+C_5 ||\f-\p||_{L^2}^{\g'},
\; \text{ where } \g'=\frac{\a}{n}\left[2/q-\g(n+2/q)\right].
$$
We finally choose $\a>0$ so large that $\g \leq \g'$
and adjust the value of the constant $C$:  this yields the desired 
estimate.
\end{proof}

Being able to control the $L^{\infty}$-norm of
$\f-\p$ by its $L^2$-norm is a powerful tool.
If for instance $\p=\f_j$, $\f$ satisfy the assumptions of proposition 3.2
-- with $\f_j$ being uniformly bounded --, and $\f_j \rightarrow \f$
in $L^1$, then $\f_j \rightarrow \f$ in 
$L^2(\o^n)$, hence $(\f_j)$ actually uniformly converges towards $\f$.
This yields the continuity of the map
$$
f \in L^p(\o^n) \mapsto \f \in {\mathcal C}^0(X),
$$
where $\f$ is the unique $\o$-psh solution to 
$(\o+dd^c \f)^n=f\o^n$, $\sup_X \f=-1$. Thus Theorem A is proved.
\vskip.1cm

We now give an application of this estimate, which requires the manifold $X$ to be 
a rational homogenous manifold, i.e. $X=G/P$ where $G$ is a complex semisimple algebraic group and $P$ a parabolic subgroup.\footnote{In particular, the cohomology class of 
$\omega$ is K\"ahler and $\omega$ itself
can be supposed to be K\"ahler without loss of generality.}.

\begin{thm}
Assume $X$ is a rational homogeneous manifold.
If $\mu=f \o^n$ is a probability measure with density $0 \leq f \in L^p(\o^n)$,
$p>1$, then the unique solution $\f \in PSH(X,\o) \in {\mathcal C}^0(X)$
to the normalized Monge-Amp\`ere equation
$$
(\o+dd^c \f)^n=\mu=f \o^n, \, \sup_X \f=-1,
$$ 
is H\"older continuous of exponent $\g>0$, for all
$\g<2/(2+nq)$, where $q=p/(p-1)$ is the conjugate exponent to $p$.
\end{thm}

\begin{proof}
Let $K$ be the maximal compact subgroup of $G$.
Then, $K$ acts transitively on $X$. Furthermore,
we may assume w.l.o.g. that $\omega$ is fixed by $K$.
One can regularize $\o$-psh functions by averaging over the Haar
measure of $K$.
This is very similar to the way one regularizes psh functions in
$\C^n$ by using convolutions with  an approximation of the
identity for the convolution product. We refer the reader to [Hu] and
the Appendix of [G] for more details.

Let $\f_{h}$ be the $\o$-psh function which is the translate
of $\f$ by an element of $K$ which is at distance $h$ from identity.
We use the notation $\f_h$ by analogy with the $\C^n$-situation,
where $\f_h(x)=\f(x+h)$.
Since $\f$ is bounded, it has gradient in $L^2$, hence
$$
||\f_h-\f||_{L^2} \leq C |h|,
$$
by using Cauchy-Schwarz inequality in a local chart.
We can thus apply proposition 3.2 to obtain that
$$
||\f_h-\f||_{L^{\infty}} \leq C' |h|^{\g},
$$
for all $\g<2/(2+nq)$. Since $\f_h(x) \simeq \f(x+h)$ in a local chart,
this precisely means that 
$\f$ is H\"older-continuous of exponent $\g$.
\end{proof}

\subsection{Regularity on the smooth locus}

\begin{thm}\label{reg}

Let $X$ be projective algebraic complex manifold, $\omega_0$ a smooth semipositive
closed $(1,1)$-form that is 
positive outside a complex subvariety $S \subset X$, and
fix $\Omega$ be a K\"ahler form on $X$.
Assume that $\omega_o^n= D \Omega^n$, where $D^{-\e}$ is in $L^1 (\Omega^n)$, and that 
$[\omega_0], [\Omega] \in NS_{\R}(X)$. 

 Let $\sigma_1, ..., \sigma_p$ (resp. $\tau_1, ..., \tau_q$) 
be holomorphic sections of some line bundle $L$ (resp $L'$) on $X$.
Fix $k\in \R_{\ge 0}$,   $l\in \R_{\ge 0}$ and $F \in {\mathcal C}^{\infty}(X,\R)$.
Assume that
$$
\int_X \frac{1}{|\tau_1|^{2l}+ \ldots + |\tau_q|^{2l}}  \Omega^n < \infty
\text{ and }
\int_X \frac{|\sigma_1|^{2k}+ \ldots + |\sigma_p|^{2k}}{|\tau_1|^{2l}+ \ldots + |\tau_q|^{2l}} e^F \Omega^n= 
\int_X \Omega^n.
$$
   
Then the unique bounded function $\f \in PSH(X,\o_0)$ such that 
$$
(\omega_0+ dd^c\f)^n = \frac{|\sigma_1|^{2k}+ \ldots + |\sigma_p|^{2k}}{|\tau_1|^{2l}+ \ldots + |\tau_q|^{2l}} 
e^{F} \Omega^n \, \text{ and } \sup_X\f=-1  
$$
is smooth outside $B=S\cup \cap_i \{\sigma_i=0\}\cup \cap_i \{\tau_i=0\} $. 
\end{thm}

\begin{rqe} This result should be compared with
[Y], Theorem 8. Yau's result is stronger in many respects
(there is no projectivity/rationality assumption and it gives a more precise 
regularity theory); 
on the other hand the conditions on the poles of the L.H.S. is less optimal than here. 

We expect the projectivity/rationality assumptions to be superfluous. 
We also  expect that a finer regularity theory might be developed for
singular  KE metrics depending
on a finer analysis of the klt singularities involved.
\end{rqe}

The rest of this subsection will be devoted to the proof of Theorem \ref{reg}. 
For the reader's convenience, we will treat two special cases
before tackling the general case
\footnote{
Notice that apart from the ${\mathcal C}^0$-estimate with degenerate L.H.S., 
the methods used here are standard and in [Y], [Ts] and [Ko]. 
Higher regularity in [TZ] is treated along similar lines given the
$L^{\infty}$-estimate the authors announce.}

\subsubsection*{Preliminary considerations}

Thanks to Lemma \ref{leb} -- here we use that $D^{-\e} \in L^1$ -- 
and Theorem \ref{mas}, for every $t\in [0,1]$ there is a 
unique continuous function
$\f_t \in PSH(X,\o_0+t\Omega)$ such that 
$$
 (\omega_o+t\Omega+ dd^c\f_t)^n = C_t\frac{|\sigma_1|^{2k}+ \ldots + 
|\sigma_p|^{2k}}{|\tau_1|^{2l}+ \ldots + |\tau_q|^{2l}} e^{F} \Omega^n
\text{ and } \sup_X \f_t=-1,
$$
where $C_t>0$ is an adequate normalisation constant and  $\|\f_t\|_{{\mathcal C}^0(X)}$
is uniformly bounded by a constant independant of $t\ge 0$.
 
We cannot use right away [Y], Theorem 8 p. 403, to ensure that $(\f_t)$ be smooth outside $B$ 
for $t>0$, since our integral condition is stronger than his. 
However we can use [Y], Thm 3, p 365 to conclude that, in case 
$\cap_i \{\tau_i=0\}=\emptyset $, $(\f_t)$ is smooth outside $B$ 
and $dd^c\f_t$ is a form whose coefficients are globally bounded on $X$. Since this does not imply ellipticity 
if $\cap_i \{\sigma_i=0\}\not=\emptyset $, this does not imply higher regularity on the whole of $X$.
 
The required uniformity in $t>0$ is not proved in [Y]. 
To deal with this case, we use a nice trick
due to H.Tsuji [Ts].

\subsubsection*{The simplest case}

First, assume $\cap_i \{\sigma_i=0\}\cup\cap_i \{\tau_i=0\}=\emptyset $. Hence the 
 family of equations under consideration  can be rewritten as: 
 
 $$
 (\omega_o+t\Omega+ dd^c\f_t)^n = C_t e^{F} \Omega^n,
 $$
$F$ being smooth. 
 
Tsuji's trick is as follows. By Kodaira's lemma,
there exists $E$ an effective Cartier divisor of $X$ such that $[\omega_o] = [\kappa_{\e}]+ \e [E]$
where $[\kappa_{\e}]$ is ample, hence we may choose a representative $\kappa_{\e}$ 
which is a K\"ahler form for every $\e >0$ small enough.  We may actually assume $E$ contains $B$ and use a family of 
$E$ such that $\cap Supp(E) =B$, by Nakamaye's theorem on base loci [Na].

Actually, despite the notation, it will NOT be necessary to let 
$\e$ decrease to $0$ \footnote{This technical device could be useful to study finer
regularity results} and we will fix once for all such an $\e>0$.

Let $\sigma\in H^0(X, \mathcal{O}_X (E))$ be the canonical section vanishing on $E$ with the appropriate
multiplicity. We can fix a smooth hermitian metric on this line bundle such that the Poincar\'e Lelong
equation holds,
$$
\omega_o= \kappa_{\e} +\e [E] - \e dd^c\log |\sigma|^2  .      
$$

The function $\phi_t:= \f_t-\e \log |\sigma|^2$ is smooth 
in $X \setminus E$ and is a classical 
solution to the PDE
$$
(\kappa_{\e}+ t\Omega + dd^c \phi_t)^n = e^{F_{\e,t}} (\kappa_{\e}+ t\Omega)^n,
$$
where $(F_{\e,t})_{1\ge t>0}$ is uniformly 
bounded in the ${\mathcal C}^{\infty}(X)$-topology of functions and $\kappa_t=\kappa_{\e}+ t\Omega$
is uniformly bounded in the ${\mathcal C}^{\infty}$-topology of K\"ahler forms on $X$.

We can use the result of the calculation in [Y], section 2. The important formula 
is (2.22) p. 351 and in a subsidiary fashion (2.21). 
In these formulae, at each point $p\in X-E$, an adequate system of normal coordinates for $\kappa_t$ is
constructed and comparing the notations here and there, we  substitute
$n$ for $m$,
$\kappa_t$ for $g_{i\bar j}$, $\kappa_t+dd^c\phi_t$ for $g_{i\bar j}'$,
$\phi_t$ for $\phi$ and  $F_{\e,t} $ for $F$.
The operator 
$\Delta$ is the Laplace operator (with the analyst's sign) of $\kappa_t$
and $\Delta'$ the Laplace operator of $\kappa_t+dd^c\phi_t$.
Also $R_{i\bar i l \bar l}=R^{t}_{i\bar i l \bar l}$ is the holomorphic bissectional 
curvature of $\kappa_t$ expressed
in the above system of normal coordinates.

Since $\kappa_t$ is uniformly bounded in the ${\mathcal C}^2$ topology of K\"ahler forms then 
certainly there is constant $C=C_{\e}$ independent of $t$ such that (2.21) holds and $C'$
also independent of $t$ such that $C'> \inf R^{t}_{i\bar i l \bar l}$. 

After these substitutions are made, (2.22) p. 351 reads:
$$
e^{C\phi_t} \Delta' (e^{-C\phi_t} (n+\Delta \phi_t) ) \ge 
\Delta (F_{\e, t} ) -n^2 C' -Cn (n+\Delta \phi_t)
+ e^{-\frac{F_{\e, t}}{n-1}} (n+\Delta \phi_t)^{\frac{n}{n-1}}
$$

We can fix  constants $C_i$ independent of $t$ such that
$$
 \Delta F_{\e,t} \ge C_1 \; \;  \; \; 
\text{ and }  \; \;  \; \;
 e^{-F_{\e, t}/n-1} \ge C_3>0.
$$
Thus setting $y=n+\Delta\phi_t$ yields
$$ 
e^{C\phi_t} \Delta' (e^{-C\phi_t} (n+\Delta \phi_t) ) \ge 
C_5 + C_6 y + e^{C_7} y^{\frac{n}{n-1}}.
$$

Now by definition
$$
e^{-C\phi_t} (n+\Delta \phi_t) = | \sigma|^{+C\e} e^{-C\f_t} (n+\Delta \f_t
+ \e \Delta \log |\sigma|^2).
$$

For each $t>0$ the functions
$\f_t$, $\e \Delta \log |\sigma|^2$ and $\Delta\f_t$ are  bounded on $X$.
Hence the positive function $e^{-C\phi_t} (n+\Delta \phi_t)$ is continuous on $X$, vanishes
on $E$ and is smooth on $X-E$. 
Its maximum is achieved at some point $p_t\not\in E$. 
It follows from the maximum principle that 
$$
0 \ge C_5 + C_6 y + e^{C_7} y^{\frac{n}{n-1}}
\text{ at point } y=y(p_t).
$$

Therefore $y\le C_8$ with a constant  independent of  $t>0$. 
Now $e^{-C\phi_t(p_t)}= |\sigma(p_t)|^{+C\e}
e^{-C\f_t(p_t)}$. Using the uniform  ${\mathcal C}^0$ estimate for $\f_t$, 
we get $0 \le  (n+ \Delta\phi_t) \le C_9 e^{+C\phi_t} $.
Since $|\f_t|$ and $\e \Delta \log |\sigma|^2$ are uniformly bounded by a 
constant independent of $t>0$, we infer
$$
(n+ \Delta\f_t) \le C_{10}  |\sigma|^{-C \e}=C_{10}  |\sigma|^{-C_{\e} \e}.
$$
This yields a $t$-independent ${\mathcal C}^0$- estimate of $dd^c\f_t$ on the compact
subsets  of 
$X-E$ \footnote{Note that $C=C_{\e}$ and that $\e C_{\e}$ 
might blow up as $\e$ goes to $0$.} . 

Standard arguments of the theory of complex Monge-Amp\`ere equations
give an interior estimate of $\f_t$ in ${\mathcal C}^{k,\alpha}_{loc}(X-E )$ 
for every $k\ge 2$, $\alpha\in]0,1[$ which is independent of $t>0$ 
(see for instance Theorem 5.1, p. 15 in [Bl2]).
 Hence the family $(\f_t)_{t>0}$ is precompact in every ${\mathcal C}^{k,\alpha}_{loc}(X-E )$. Its 
 cluster values are cluster values in ${\mathcal C}^0(X-E)$ hence they are all equal to $\f|_{X-E}$. 
 This implies $\f\in {\mathcal C}^{k,\alpha}_{loc}(X-E)$, 
hence that $\f\in {\mathcal C}^{\infty}(X-E)$.

\subsubsection*{Case where $\cap_i \{\tau_i=0\}=\emptyset $}
\footnote{It suffices to consider this case for constructing singular KE metrics on 
algebraic varieties with canonical singularities}
We study here the equation
$$
 (\omega_o+t\Omega+ dd^c\f_t)^n = C_t(|\sigma_1|^{2k}+ \ldots + |\sigma_p|^{2k}) e^{F} \Omega^n
$$
 
The first few steps of the preceding argument can be repeated without changes. 
Next we apply formula (2.22) in [Y] as earlier, except that we set 
$F= F_{\e,t} + \log ||s||^{(2k)}$, where  
$||s||^{(2k)}:=|\sigma_1|^{2k}+ \ldots + |\sigma_p|^{2k}$.  This yields
\begin{eqnarray*}
e^{C\phi_t} \Delta' (e^{-C\phi_t} (n+\Delta \phi_t) ) 
&\ge & \Delta F_{\e, t}  + \Delta \log ||s||^{(2k)}  
-n^2 C' -Cn (n+\Delta \phi_t) \\
 & & + \left( \frac{e^{-F_{\e, t}}}{||s||^{(2k)}} \right)^{1/(n-1)}
(n+\Delta \phi_t)^{\frac{n}{n-1}}
\end{eqnarray*}

We recall the two preceding inequalities and observe two new ones that are available:
\begin{eqnarray*}
\Delta F_{\e,t} \ge C_1 \; & \text{ and }& \; 
 e^{-F_{\e, t}/n-1} \ge C_3>0; \\
\Delta  \log ||s||^{(2k)} \ge C_2  
&\text{ and }&
C_4 \ge ||s||^{(2k)}.
\end{eqnarray*}
 
Setting as earlier $y=n+\Delta\phi_t$, we get 
$$ 
e^{C\phi_t} \Delta' (e^{-C\phi_t} (n+\Delta \phi_t) ) \ge 
C_5 + C_6 y + e^{C_7} y^{\frac{n}{n-1}}.
$$

After this point, the proof is entirely the same as before. 

\begin{rqe}
In order to carry out the second order a priori estimate, one needs 
information that only depend on $\sup_X F$ and $\inf_X \Delta F$.
This is pointed out in [Y], p. 351, and it is the basis 
for the proof of [Y], Thm 3. 
\end{rqe}

\subsection{Formal reduction of the general case to the second case using 
smooth orbifolds}

In order to carry out the present argument, which in essence is just 
a change of variables $z \rightsquigarrow \zeta =z^{1/m}$, we 
need to use analysis on certain smooth orbifolds. We will not give
complete definitions since they are in the
recent reference [BGK], section 2 pp. 560-564, see also [MO] and the references therein. 

Let $(X,\Delta)$ be a smooth orbifold pair. By this we mean 
that we have the prime decomposition $\Delta = \sum_i (1-m_i^{-1}) E_i$
where $m_i\in \N^*$ is an integer. We assume that
$supp(\Delta)$ is a simple normal crossing divisor.
 Then, a classical construction
surveyed in  [BGK] enables to construct an orbifold
$[X,\Delta]$ with a $1$-morphism of orbifolds
$c:[X,\Delta] \to X$ with the following properties:

\begin{itemize}
\item $c$ is the reduction to the coarse moduli space of $[X,\Delta]$.
\item $c_{X-Supp(\Delta)}: U:= [X,\Delta]\times_X  ( X-Supp(\Delta) )\to X-Supp(\Delta)$ is an isomorphism. 
Hence $U$ is an open suborbifold of $[X,\Delta]$ which is 
an old-fashioned manifold).
\item For every open polydisk $\mathbb{D}\subset X$ with local coordinates
$z_1,\ldots, z_n$  such that $Supp(\Delta)= \{\prod_{j=1}^p z_j =0\}$
$[X,\Delta]\times_X \mathbb{D}= [\mathbb{D}'/ G_{loc}]$.

In this formula, the local 
isotropy group is
$G_{loc}=\prod_{j=1}^p \Z / m_j \Z$, $m_j$ is the integer 
multiplicity of the divisor $E_{i_j}$ such that $E_{i_j}\cap \mathbb D
=\{ z_j=0\}$, $G_{loc}$ acts on the polydisk $\mathbb {D}'$ by 
$(\zeta_1,...,\zeta_p).(z'_1, ..., z'_n)= (\zeta_1 z'_1, ..., \zeta_p z'_p,
z_{p+1}', ...)$ \footnote{The usual isomorphism of $\Z/ m\Z$ with the group 
of $m$-th root of unity is used.}. 

The orbifold $1$-morphism $[\mathbb{D} '/ G_{loc}]\to \mathbb{D}$ is induced by 
$\varkappa_{loc}: (z'_1, .., z'_n)\mapsto ((z'_1)^{m_1}, ...)$. 

\item For sufficiently divisible $s$,  $c_* \mathcal{O}_{[X,\Delta]} (s K_{[X,\Delta]})=
\mathcal{O}_X (s(K_X+\Delta))$.
\end{itemize}

It is possible to define all the basic concepts of K\"ahler geometry on orbifolds
such as smooth functions, K\"ahler metrics, etc... The principle is to think of $\varkappa_{loc}^{-1}$
as a (multivalued) smooth coordinate chart.

A continuous function on $[X,\Delta]$ is a continuous function on $X$. A Radon measure on $[X,\Delta]$
is a Radon measure on $X$.

A smooth function $f$ on $[X,\Delta]$ is a continuous function on $X$ such that for every local chart 
$\varkappa_{loc}^* f$ is smooth. In particular $f$ is H\"older continuous. 
 
A K\"ahler metric $\Omega_{[X,\Delta]}$ on $[X,\Delta]$
is a K\"ahler metric $\Omega_{X-Supp(\Delta)}$ on $X-Supp(\Delta)$ with 
the property that $\varkappa_{loc}^*\Omega$ extends to a smooth K\"ahler metric on $\mathbb{D}'$. 
In particular, it also extends as a closed K\"ahler current on $X$ with H\"older potentials. 

The pull back of a K\"ahler form on $X$ to $[X,\Delta]$ is a 
semipositive closed $(1,1)$- form that is actually  cohomologous to
a K\"ahler class\footnote{Here no reference can be given.  But it is easy to extend 
the gluing methods for K\"ahler forms 
developed in [Dem 3] and 
[Pa]  to orbifolds. Hence [DP] extends  to K\"ahler orbifolds. }. 

Observe that $\varkappa ^*_{loc} dz_l=  m_{i_l} (z'_l)^{m_{i_l}-1} dz'_l$, hence a smooth 
volume form  on  $[X,\Delta]$ can be interpreted as a volume form $v$ on $X-Supp(\Delta)$
such that
$$
v \text{ is comparable to }
\frac{ {\prod_{l=1}^n (\sqrt{-1} dz_l\wedge d\bar z_l)}}{\prod_{l=1}^p |z_l|^{2(1-1/m_{i_l})}}.
$$ 

In case the pair $(X,l^{-1} (\tau_1))$ is an orbifold pair, the equation 
\begin{equation} \label{eqst}
(\omega_o+ dd^c\f_t)^n = C_t\frac{|\sigma_1|^{2k}+ \ldots + |\sigma_p|^{2k}}{|\tau_1|^{2l}} e^{F} \Omega_X^n 
\end{equation}
can be interpreted on $[X,\Delta]$ as an equation of the form
$$
(c^*\omega_o+ dd^c\f_t)^n = C_t (|\sigma_1|^{2k}+ \ldots + |\sigma_p|^{2k})e^{F} \Omega_{[X,\Delta]}^n.
$$
 
 The method used to analyze the case where $\cap_i \{\tau_i=0\}=\emptyset $ extends with 
 almost no changes to the orbifold case. Hence the unique continuous solution of equation
 (\ref{eqst}) is smooth outside its singular locus if $(X,\frac{1}{l} (\tau_1))$ is an orbifold pair. 
 
 Under the more general hypothesis that $\int_X |\tau_1|^{-2l}<\infty$ and $\frac{1}{l}(\tau_1)$ is  a 
 divisor with simple normal crossings, then 
 we can construct an orbifold pair $(X,\Delta)$ with $0\le \frac{1}{l}(\tau_1) \le \Delta$ and 
 we are back to the previous case. 
 
 For the most general case, consider the ideal $\mathcal{I}$ generated by the $t_i$
 and fix $\mu: X'\to X$ a log resolution of $(X,\mathcal{I})$. Then we are back to the 
 previous case, with an equation
 on $X'$ . 
This ends the proof of Theorem 3.6. 
\vskip.2cm

In certain rare circumstances,
 there is a finite smooth covering $Y\to X$ such that $Y/G =X$ and $[Y/G]=[X,\Delta]$ and
the argument we use here reduces to a $G$-equivariant argument on $Y$.

\begin{rqes}
If we start with $\omega_o$ K\"ahler, and the log-resolution is non trivial, 
$\mu^*\omega_o$ is not K\"ahler anymore. 
  
This method that dates back to [Ko] can be used to prove 
a variant of [Y], Theorem 7 p. 399 where the divisor of $\sigma_2$ is 
a simple normal crossing divisor, under the sole assumption that 
$\int_M |\sigma_2|^{-2k_2} <\infty$.

Now, it could not have been used to prove Theorem 8 p. 403 in 1978 since 
 log-resolutions force the use of
 Monge-Amp\`ere equations with degenerate L.H.S, for which the ${\mathcal C}^0$-estimate
 proved here
 was not available then.
\end{rqes}

\section{More Monge-Amp\`ere equations}

As we aim at constructing singular K\"ahler-Einstein metrics, 
it is important to consider Monge-Amp\`ere equations of the following type,
$$
(\o+dd^c \f)^n=e^{t \f} \mu,
$$
where $\mu$ is a probability measure which satisfies condition
${\mathcal H}(\a,A,\o)$ (see definition 1.3), and $t$ is a real parameter.
The case $t=0$, treated in Theorem 2.1, will correspond to Ricci-flat metrics
(see section 6). We focus here on case 
$t>0$.

\begin{thm}\label{mma}
Let $\mu$ be a probability measure which satisfies condition
${\mathcal H}(\a,A,\o)$ and fix $t>0$.
There exists a unique function $\f_t \in PSH(X,\o) \cap L^{\infty}(X)$
such that
$$
(\o+dd^c \f_t)^n=e^{t \f_t} \mu.
$$
Furthermore $\f_t$ is continuous provided $(X,[\omega])$ enjoys the continuous approximation property.
\end{thm}

\begin{proof}
The uniqueness easily follows from the comparison principle as
we explain in proposition 4.3 below.
We are going to prove the existence by a fixed point method.

Fix $\p \in {\mathcal E}^1(X,\o)$ such that $\int_X \p d\mu=0$,
and let us consider the equation
$$
MA(\p) \hskip2cm (\o+dd^c \f)^n=e^{t\p-c_{\p}} \mu,
$$
where the constant $c_{\p}:=\log [\int_X e^{t\p} d\mu ]$
is chosen so that 
$$
1=\int_X (\o+dd^c \f)^n=e^{-c_{\p}} \int_X e^{t\p} d\mu.
$$
Observe that $\mu_{\p}:=e^{t\p-c_{\p}} \mu$ satisfies condition
${\mathcal H}(\a,A_{\p},\o)$,
where $A_{\p}=\exp( t \sup_X \p -c_{\p})$.
It follows therefore from Theorem 2.1 that there exists a
unique bounded function $\f \in PSH(X,\o)$ solution
to $MA(\p)$ and normalized by $\int_X \f d\mu=0$.
We use here this {\it linear} normalization rather than the 
non-linear $\sup$-normalization: they are comparable
thanks to proposition 2.7 in [GZ 1], which shows that
$$
-M_{\mu} \leq \int_X u d\mu -\sup_X u \leq 0,
$$
for all functions $u \in PSH(X,\o)$ and for
some uniform constant $M_{\mu}>0$.
Since $\int_X \p d\mu=0$, we infer
\begin{equation}
0 \leq {\mathcal E}_{\o}(\f):=\int_X |\f| \o_{\f}^n=e^{-c_{\p}} \int_X |\f| e^{t\p} d\mu
\leq 2 M_{\mu} e^{t M_{\mu}},
\end{equation}
by observing that $c_{\p} \geq 0$ since $t \geq 0$, and
$$
\int_X |\f| d\mu \leq \int_X |\f - \sup_X \f| d\mu +\sup_X \f \leq 2 M_{\mu},
$$
since $\int_X \f d\mu=0$.

The important fact here is that the energy ${\mathcal E}_{\o}(\f)$ of $\f$ is bounded
from above by a constant $M_0:=2M_{\mu}e^{t M_{\mu}}$ which is independent of $\p$.
We have thus defined an operator
$$
T:\p \in {\mathcal C}_M \mapsto \f \in {\mathcal C}_{M_0}
$$
which associates to $\p \in {\mathcal C}_M$ the unique solution
$\f \in {\mathcal C}_{M_0}$ to $MA(\p)$, where
$$
{\mathcal C}_{M}:=\left\{ \p \in {\mathcal E}^1(X,\o) \, / \, 
\int_X \p d\mu=0 \text{ and } {\mathcal E}_{\o}(\p) \leq M \right\}.
$$
It follows from [GZ 2] that
${\mathcal E}^1(X,\o)$ is convex. So is the subset of functions 
$\p \in {\mathcal E}^1(X,\o)$
such that $\int_X \p d\mu=0$.
The set ${\mathcal C}_M$ is not convex, but it is relatively
compact in $L^1(X)$ and its closed convex hull
$\hat{{\mathcal C}}_M$ is contained in
${\mathcal C}_{\kappa_n M}$ for some uniform constant 
$\kappa_n$ which only depends on the dimension of $X$:
this follows from easy computations (see 
proposition 2.10 in [GZ 2]).
Therefore $T$ maps the compact convex set $\hat{{\mathcal C}}_M$
into itself if $M$ is large enough.

We claim that $T$ is continuous.
Let $(\p_j) \in {\mathcal C}_M^{\N}$ be a sequence of functions which converges
in $L^1(X)$ towards $\p \in {\mathcal C}_M$.
We need to show that $\f_j:=T(\p_j)$ converges in $L^1(X)$ towards $T(\p)$.
Since the set $\{ u \in PSH(X,\o) \, / $
$\, \int_X u d\mu=0 \}$ is relatively 
compact in $L^1(X)$ (see proposition 2.7 in [GZ 1]), we can assume -- relabelling
if necessary -- that $(\f_j)$ converges in $L^1(X)$ towards a function $\f \in PSH(X,\o)$.
We show in lemma 4.2 below that $(\f_j)$ converges in $L^1(\mu)$ towards $\f$.
In particular $\int_X \f d\mu=0$ and,
passing to a subsequence if necessary, we can assume that
$e^{t \p_j(x)} \rightarrow e^{t \p(x)}$
for $\mu$ almost every point $x$.
Set
$$
\hat{\f}_j:=\left( \sup_{l \geq j} \f_l \right)^* \; \; 
\text{ and } \; \; 
\check{\p}_j:=\inf_{l \geq j} \p_l.
$$
Observe that $(\hat{\f}_j)$ decreases towards $\f$, while
$(e^{t \check{\p}_j})$ increases towards $e^{t \p}$ at $\mu$ almost 
every point. The energy of $\hat{\f}_j$ is controlled by that of $\f_j$
since $\hat{\f}_j \geq \f_j$ (see lemma 2.3 in [GZ 2]), and
${\mathcal E}_{\o}(\f_j) \leq M_0$ by (5), therefore
$\f \in {\mathcal E}^1(X,\o)$ and
$(\o+dd^c \hat{\f}_j)^n \rightarrow (\o+dd^c \f)^n$.
It follows from an inequality of J.-P.Demailly [Dem 1] that
$$
(\o+dd^c \hat{\f}_j)^n \geq e^{t\check{\p}_j -\hat{c}_j} \mu,
$$
where $\hat{c}_j:= \sup_{l \geq j} c_{\p_l}$.
Observe that $\hat{c}_j \rightarrow c_{\p}$, thus
$$
(\o+dd^c \f)^n \geq e^{t \p-c_{\p}} \mu.
$$
Since these are two probability measures, there is actually equality
hence $\f=T(\p)$: this shows that $T$ is continuous.

We can now invoke Schauder fixed point theorem, which yields a fixed
point $\f=T(\f), \f \in {\mathcal C}_M$.
The function $\f$ is automatically bounded (by Theorem 2.1, since
$e^{t\f-c_{\f}} \mu$ satisfies ${\mathcal H}(\a,A',\o))$, hence
$\Phi:=\f-t^{-1} c_{\f}$ is the solution we were looking for.
\end{proof}

\begin{lem}
The functions $\f_j=T(\p_j)$ (respectively $e^{t \p_j}$) converge in $L^1(\mu)$ towards $\f$
(respectively $e^{t \p}$).
\end{lem}

\begin{proof}
We first show that $(\f_j)$ converges to $\f$ in $L^1(\mu)$.
Observe that the sequence $(\f_j)$ is uniformly bounded: this follows
from Theorem 2.1 since $(\o+dd^c \f_j)^n$ satisfies ${\mathcal H}(\a,A_j,\o)$,
where $A_j=e^{t \sup_X \p_j-c_{\p_j}} A \leq e^{t M_{\mu}} A$ is bounded
from above.
It follows then from standard arguments that 
$\int_X \f_j d\mu \rightarrow \int_X \f d\mu$ (see e.g. the proof of lemma 5.2 in [Ce]).

Fix $\e>0$ and let $G$ be an open set of $X$ such that $\f$ is continuous on
$X \setminus G$ and $Cap_{\o}(G) \leq \e$ (see corollary 3.8 in [GZ 1]).
By Hartogs' lemma, $\f_j \leq \f+\e$ on the compact set $X \setminus G$,
if $j \geq j_{\e}$. Observe that
$$
\int_{X \setminus G} |\f-\f_j| d\mu \leq 2 \e+\int_{X \setminus G} [\f-\f_j] d\mu
\leq 3 \e,
$$
if $j \geq j'_{\e}$. O the other hand since $\mu$ satisfies ${\mathcal H}(\a,A,\o)$,
we get
$$
\int_G |\f-\f_j| d\mu \leq 2M \mu(G) \leq 2M A \e^{1+\a},
$$
where $M=\sup_j ||\f_j||_{L^{\infty}(X)}$. This shows that 
$||\f-\f_j||_{L^1(\mu)} \rightarrow 0$.

The proof for $(e^{t \p_j})$ is similar: it suffices to note that
the functions $u_j:=e^{t\p_j-t \sup_X \p_j}$ 
are $\o$-psh and uniformly bounded.
One can then apply the rest of the argument.
\end{proof}

\begin{pro}
Let $\f,\p \in {\mathcal E}^1(X,\o)$ and $t>0$ be such that
$$
(\o+dd^c \f)^n e^{-t\f} =(\o+dd^c \p)^n e^{-t\p}.
$$
Then $\f \equiv \p$.
\end{pro}

\begin{proof}
It follows from the comparison principle (see [GZ 2]) that
$$
\int_{(\f<\p)} (\o+dd^c \p)^n \leq \int_{(\f<\p)} (\o+dd^c \f)^n 
=\int_{(\f<\p)} e^{t(\f-\p)}(\o+dd^c \p)^n .
$$
Since $e^{t(\f-\p)}<1$ on $(\f<\p)$, we infer $\f \geq \p$ for
$\nu$ almost every point, where $\nu=(\o+dd^c \p)^n$.
Reversing the roles of $\f,\p$ yields $\f = \p$ for $\nu$
almost every point. Therefore $(\o+dd^c \f)^n=(\o+dd^c \p)^n$,
hence $\f-\p=c$ is constant by Theorem 3.4 in [GZ 2].
Finally $c=0$ since $e^{tc}=1$ and $t>0$.
\end{proof}
\vskip.2cm

When $K_X$ is nef and big, H.Tsuji constructed in [Ts] -- using K\"ahler-Ricci
flow techniques -- a function $\p \in PSH(X,\o)$ such that
$$
\int_X e^{t \p} d\mu=1, \; 
\p \in {\mathcal C}^{\infty}(X \setminus E) \;
\text{ and } (\o+dd^c \p)^n=e^{t \p} \mu
\text{ in } X \setminus E,
$$
where $E$ is the
exceptionnal locus of the map associated to
the base point free linear system $|NK_X|$, $N\in\N$ big enough and
the current $T_{KE}=\o+dd^c \p$  defines a K\"ahler-Einstein metric. This function coincides with our solution thanks to the following
unicity result.

\begin{pro}
Let $\mu$ be a probability measure and $t>0$. Let $\f,\p \in PSH(X,\o)$ be such that
$\int_X e^{t \f} d\mu=\int_X e^{t \p} d\mu=1$.
Assume $\f \in {\mathcal E}^1(X,\o)$ is a global solution
to the complex Monge-Amp\`ere equation $(\o+dd^c \f)^n=e^{t \f} \mu$,
while $\p \in {\mathcal C}^0(X \setminus E)$
satisfies $(\o+dd^c \p)^n=e^{t \p} \mu$ only in $X \setminus E$.

Then $\p \in {\mathcal E}^1(X,\o)$ and $\p \equiv \f$. 
\end{pro}

\begin{proof}
Set $\p_j:=\max(\p,-j) \in PSH(X,\o) \cap L^{\infty}(X)$.
Observe that the probability measures $(\o+dd^c \p_j)^n$ converge in $X \setminus E$
towards the measure $\nu=e^{t \p} \mu$. Since 
$\nu(X)=\nu(X \setminus E)=1$, it follows that 
$(\o+dd^c \psi_j)^n$ converges to $\nu$ on all of $X$.

Fix $\e>0$ and set $v_{\e}:=(\p+\e v)/(1+\e) \in PSH(X,\o)$, where 
$v \in PSH(X,\o)$, $v \leq 0$, is such that $e^v $ is continuous and $(v=-\infty)=E$.
It follows from lemma 2.2 that
for all $s>0$,
$$
Cap_{\o}(\p_j<-s-1) \leq \int_{(\p_j<-s)} (\o+dd^c \p_j)^n
\leq \int_{(v_{\e} \leq -s/(1+\e))} (\o+dd^c \p_j)^n.
$$
Observe that $e^{v_{\e}}$ is continuous on $X$, hence the sublevel sets
$(v_{\e} \leq c)$ are compact.
We infer, letting $j \rightarrow +\infty$,
$$
Cap_{\o}(\p<-s-1) \leq \int_{(v_{\e} \leq -s/(1+\e))} e^{t \p} d\mu.
$$
Letting $\e$ go to zero and using that $\mu(X)=1$ yields
$$
Cap_{\o}(\p<-s-1) \leq \int_{(\p<-s)} e^{t \p} d\mu \leq e^{-s}.
$$
Therefore the capacity of the sublevel sets of $\p$ decreases fast
as $s \rightarrow +\infty$, hence by
lemma 5.1 in [GZ 2] we get $\p \in {\mathcal E}^1(X,\o)$.
Since $e^{-t\f}(\o+dd^c \f)^n \equiv e^{-t\p}(\o+dd^c \p)^n$, it follows 
from proposition 4.3 that $\f \equiv \p$.
\end{proof}

\vskip.2cm

\begin{thm}\label{reg2}

Let $X$ be projective algebraic complex manifold, $\omega_0$ a smooth semipositive closed $(1,1)$-form that is 
K\"ahler outside a complex subvariety $S \subset X$, and
fix $\Omega$ be a K\"ahler form on $X$.
Assume that $\omega_o^n= D \Omega^n$, where $D^{-\e}$ is in $L^1 (\Omega^n)$, and that 
$[\omega_0], [\Omega] \in NS_{\R}(X)$. 

 Let $\sigma_1, ..., \sigma_p$ (resp. $\tau_1, ..., \tau_q$) 
be holomorphic sections of some line bundle $L$ (resp $L'$) on $X$.
Fix $k\in \R_{\ge 0}$,   $l\in \R_{\ge 0}$ and $F \in {\mathcal C}^{\infty}(X,\R)$.
Assume that
$$
\int_X \frac{1}{|\tau_1|^{2l}+ \ldots + |\tau_q|^{2l}}  \Omega^n < \infty
\text{ and }
\int_X \frac{|\sigma_1|^{2k}+ \ldots + |\sigma_p|^{2k}}{|\tau_1|^{2l}+ \ldots + |\tau_q|^{2l}} e^F \Omega^n= 
\int_X \Omega^n.
$$
   
For each $t>0$, the unique function $\f \in PSH(X,\o_0) \cap L^{\infty}(X)$ such that 
$$
(\omega_0+ dd^c\f)^n = \frac{|\sigma_1|^{2k}+ \ldots + |\sigma_p|^{2k}}{|\tau_1|^{2l}+ \ldots + |\tau_q|^{2l}} 
e^{F+t \f} \, \Omega^n 
$$
is smooth outside $B=S\cup \cap_i \{\sigma_i=0\}\cup \cap_i \{\tau_i=0\} $. 

\end{thm}

\begin{proof}
The proof of Theorem \ref{reg} applies here almost verbatim.  
\end{proof}

\begin{rqe}
We will apply Theorem 4.1 in section 6 to construct singular 
K\"ahler-Einstein metrics
on manifolds of general type.
This will follow from the resolution of $(\o+dd^c \f)^n=e^{t \f} \mu$ for
large enough values of $t>0$.
The Monge-Amp\`ere equations
$$
(\o+dd^c \f)^n=e^{-t \f} \mu, \; \; t>0,
$$
can also be solved with a similar method, but only for small values of $t<t_X$.
The critical exponent $t_X$ depends on the manifold $X$, and 
may be too small to produce K\"ahler-Einstein metrics when $c_1(X)>0$:
even smooth
manifolds of positive scalar curvature do not necessarily admit K\"ahler-Einstein
metrics (see [T]). Since technical details are much more involved in this case,
we postpone this study to a forthcoming article.
\end{rqe}

\section{Singularities in Mori theory}

The singular locus of the normal complex space of pure dimension
$n$ is a codimension $\ge 2$ analytic subvariety denoted by $V^{sing}$. Let
$V^{reg}= V-V^{sing}$ and $j:V^{reg}\to V$ be the natural open immersion.

\subsection{Log terminal singularities}

Since this material may not be familiar to complex analysts
 or differential geometers, we briefly recall some basic facts
  on some of the singularities encountered in the 
Minimal Model Program (MMP for short). See [KM]
for a detailled account in the algebraic case, the analytic theory being also surveyed
 there in less detail.

The sheaf of holomorphic functions $\mathcal O_V$ is the subsheaf of the sheaf of continuous functions on $V$
consisting of the functions whose restriction to $V^{reg}$ is holomorphic. Actually, by Hartogs' theorem, any holomorphic function
on $V^{reg}$ extends to $V$, which means that $j_*\Oxr=\Ox$.

Every meromorphic n-form $\alpha$ on $V^{reg}$ extends to  $V$,
i.e. let $\pi:X\to V$ be a resolution of singularities of $V$, 
then the meromorphic $n$-form $\pi^* \alpha $ defined on $\pi^{-1} V^{reg}$ 
extends to a meromorphic $n$-form on $X$.
Let $\omega_{V^{reg}}$ be the canonical sheaf of the smooth variety $V^{reg}$. 
The sheaf
$\omega_V=j_*\omega_{V^{reg}}$ is a coherent analytic sheaf on $V$. 

More generally every meromorphic pluricanonical form on $V^{reg}$ extends to $V$ and
$\omega^{[q]}_V= j_*\omega_{V^{reg}}^q$, $q>0$ is a coherent analytic sheaf on $V$.

\begin{defi}\label{gor}
Say $V$ is {\bf 1-Gorenstein} iff one of the following equivalent conditions holds:
\begin{enumerate}
\item Every $x\in V$ has an open neighborhood $U$ such that $U^{reg}$ carries a holomorphic $n$-form with an empty zero divisor.
\item  $\omega_V$ is a rank one locally free sheaf. 
\item  Every $x\in V$ has an open neighborhood $U$ such that $\omega_{U^{reg}}$ is isomorphic to $\Oxr|_U$.
\end{enumerate}
A local section of $\omega_V$ defining a holomorphic $n$-form without zeroes on $V^{reg}$ will be called a 
local generator of $\omega_V$. If furthermore  $V$ is Cohen-Macaulay, $V$ is said to be Gorenstein. 

Say $V$ is {\bf $\Q$-Gorenstein} iff one of the following equivalent conditions 
is satisfied:

\begin{enumerate}
\item Every $x\in V$ has an open neighborhood $U$ such that $U^{reg}$ carries a holomorphic pluricanonical  form with an empty zero divisor.
\item  For every $x\in V$, there exists $ N_x \in\N$ and an open neighborhood 
$U$ of $x$ such that $\omega^{[N_x]}_U$ is a rank one locally free sheaf.
\item  For every $x\in V$ there is $N_x\in \N$ and   an open neighborhood $U$ of $x$ such that $\omega^{N_x}_{U^{reg}}$ is isomorphic to $\Oxr|_U$.
\end{enumerate}

A local section of $\omega^{[N]}_V$ defining a holomorphic pluricanonical form without zeroes on $V^{reg}$ will be called a 
local generator of $\omega^{[N]}_V$. 

For every $x\in V$, the smallest $N_x$ fulfilling condition 3 near $x$ is called the local index of $V$ at $x$. The l.c.m. of all local indices, if finite, is called the index of $V$. 
\end{defi}

\begin{defi} \label{can}
Say $V$ has only {\bf canonical singularities} iff $V$ is $\Q$-Gorenstein, of finite index $N$ and one of the following equivalent conditions is fulfilled:
\begin{enumerate}
\item Let $\pi:X\to V$ be a resolution. Let $\alpha$ be a local generator of $\omega_V^{[N]}$.
The meromorphic pluricanonical form $\pi^*\alpha$ is holomorphic.
\item Let $\pi:X\to V$ be a resolution. For every $m\in \N, \  \pi_* \omega_{X}^{[Nm]}=\omega_V^{[Nm]}$.
\item (Assuming $V$ is an algebraic variety) 
Let $\pi:X\to V$ be a resolution. 
Then $K_{X}\cong \pi^* K_V + \displaystyle\sum  a_E E$ with $a_E\ge 0$ 
where $\cong$ means numerical equivalence of $\Q$-Cartier divisors
and the sum runs over the exceptional divisors of $\pi$.
\end{enumerate}
\end{defi}

Observe that it is enough to check the first two conditions for {\it some} resolution. 
In the third condition $Na_E$ is the order of vanishing of $\pi^*\alpha$ along the divisor $E$.

\begin{defi}
Say $V$ has only {\bf log-terminal singularities} iff $V$ is  $\Q$-Gorenstein, 
of finite index  $N$ and the following holds:
let $\pi:X\to V$ be a log-resolution and let $\alpha$ be a local generator of $\omega_V^{[N]}$:
then the pole along any component $E$ of $exc(\pi)$ of the meromorphic $N$-canonical form  
$\pi^*\alpha$ on $X$ is of order $\le N-1$.

When $V$ is algebraic an equivalent formulation is: let  $\pi:X\to V$ be a  log-resolution.
Then $K_{X}\cong \pi^* K_V + \displaystyle\sum_{E } a_E E$ with $a_E>-1$.
\end{defi}

The importance of the class of canonical singularities
\footnote{On the other hand, the class of log-terminal singlarities of varieties is less important. Indeed let $X$ be a complex projective normal  variety with log terminal singularities, then there is a Deligne Mumford stack $\mathcal{X}\to X$ which is etale in codimension one and has only Gorenstein canonical singularities. At the expense of working with this canonical cover, one could avoid  the consideration of log terminal singularities for varieties. } comes from a 
theorem due to M. Reid [R 1] (see also [Deb], p. 174):

\begin{thm}
 \label{Reid}
Let $X$ be a projective algebraic manifold of general type
 whose canonical ring $R={\displaystyle\oplus_{n\in \N}} H^0(X,\omega_X^n)$ is of finite type. 
Then the canonical model of $X$,  $X_{can}:=Proj(R)$ 
has only canonical singularities. If $N=Index( X_{can})$ then 
$\omega^{[N]}_{X_{can}} $ is ample.
\end{thm}

The finiteness of the canonical ring for varieties of general type is known in dimension 3 [Ka].
In higher dimension, Y.Kawamata has proved that it is a consequence of the existence of minimal
models. 
$X_{can}$ is a uniquely defined singular birational model of $X$. The minimal models of $X$
in the sense of the MMP are crepant terminalizations of $X_{can}$ and do not enjoy the
above unicity since they may be related by non trivial flops.

\begin{exas} Let  $S$ be a normal algebraic surface. The following are equivalent:

\begin{enumerate}
\item $S$ has only  canonical singularities.
\item $S$ is locally analytically isomorphic to
 $X=\C^2 /G$, $G\subset SL_2(\C)$ a finite subgroup.
\item The exceptional divisors of the minimal resolution $\pi_{min}$ of $S$,  
 have simple normal crossings, their
 components are  (-2) smooth rational curves, their  incidence graphs are of type 
  A-D-E (Du Val singularities).
\end{enumerate}

The log terminal surface singularities are precisely the singularities of the form 
$X=\C^2 /G$, $G\subset GL_2(\C)$ a finite subgroup.
\end{exas}

\begin{exas} 
In higher dimension, quotient singularities are still log terminal. 
Fix $n>0$ and
let $H\subset \C{\mathbb P}^{n+1}$ be a smooth degree $d$ hypersurface. The affine cone over 
$H$ has only canonical singularities iff $d\le n+1$.

In particular, the ordinary double point $x^2+y^2+z^2+t^2=0$ has only canonical singularities
but it is not a quotient singularity. 

The hypersurface singularities of type $A-D-E$ are canonical. 
\end{exas}

\subsection{Normal K\"ahler spaces}

\subsubsection*{Plurisubharmonic functions}

Let $V$ be a normal analytic space of pure dimension  $n$.
A plurisubharmonic (psh) function $\f$ on $V$  is an upper semicontinuous function on $V$
with values in $\R\cup\{-\infty\}$, which is not locally $-\infty$, and extends to a 
psh function in some local embedding $V\to \C^N$. 
The function $\f$ is strongly psh  
(resp. ${\mathcal C}^0$, resp. ${\mathcal C}^{\infty}$) iff 
it extends to a strongly psh function (resp. ${\mathcal C}^0$, resp. 
${\mathcal C}^{\infty}$) in some local embedding.
A continuous function is psh iff  its restriction to  $V^{reg}$ is so [FN].
 A bounded psh function on $V^{reg}$
extends to $V$. 

A pluriharmonic function on $V$  is a real valued continuous function on $V$
 $f$ on $V$  such that one of the following equivalent conditions holds:
\begin{itemize}
\item $f$ is locally the real part of a holomorphic function.
\item Given a local embedding $V\to \C^N$, $f$ extends locally to a pluriharmonic function on $\C^N$.
\item $f|_{V^{reg}}$ is pluriharmonic.
\end{itemize}

\subsubsection*{Semi-K\"ahler currents}

\begin{defi}
A semi-K\"ahler, resp. K\"ahler, resp. smooth K\"ahler, potential on $V$ is a 
family $(U_i,\f_i)_{i\in I}$ where $(U_i)$
is an open covering of $V$ and $\f_i$ a  psh function, 
resp. a strongly psh function,
resp. a ${\mathcal C}^{\infty}$-smooth strongly psh function,
on  $U_i$ such that $\f_i-\f_j$ is pluriharmonic
on $U_ i \cap U_j$.
\end{defi}

Define an equivalence relation on semi-k\"ahler potentials 
requiring that $(U_i,\f_i)\sim(V_j,\psi_j)$
iff $\f_i-\psi_j$ is pluriharmonic on $U_i\cap V_j$.

\begin{defi}
A smooth K\"ahler metric $\Omega$ on  $V$ is a 
$\sim$-equivalence class of smooth K\"ahler potentials.
A semi-K\"ahler (resp. K\"ahler) 
current on $V$ is a $\sim$-equivalence class of semi-K\"ahler (resp. K\"ahler) potentials.

A semi-K\"ahler current $\Omega$  is said to have 
$L_{loc}^{\infty}$ (resp. ${\mathcal C}^0$, resp. H\"older continuous) potentials iff, given a potential $ (U_i,\f_i)_{i\in I}$ for $\Omega$, each 
$\f_i$  is $L_{loc}^{\infty}$ (resp. ${\mathcal C}^0$, resp. H\"older continuous).
\end{defi}

We will on occasion drop the requirement that the local potentials of $\Omega$ are psh, replacing it by the
requirement that they are locally the sum of a smooth and a psh function. 
The current $\Omega$ will  then be called
a quasi positive closed current on $V$.

 If it has locally bounded potentials, $\Omega$ is fully determined by the  closed $(1,1)$ form  $\Omega_{reg}$ on $V_{reg}$ defined on $U_i$ by $\Omega_{reg}=dd^c\f_i$.

Let $\Omega$ be a smooth K\"ahler metric on $V$ with  K\"ahler potential
 $(U_i,\f_i)$. An upper semi-continuous function $\f: X\to \R\cup{-\infty}$ is said to be 
$\Omega$-psh iff $\forall i$ $\f_i+\f$ is psh on $U_i$.  The semi-K\"ahler current
whose potential is $(U_i,\f+\f_i)$ is denoted by $\Omega+dd^c\f$.

\begin{exa}
Let $V= \C^2/{\pm 1}$. Let $(x,y)$ be the usual affine coordinates on  $\C^2$, $(u,v,w)$
those on  $\C^3$.
The formulas $u=x^2, \ v=y^2, \ w=xy$ realize $V$ as the closed subscheme of $\C^3$
whose equation is $uv-w^2=0$.  We have two \lq natural\rq \ 
K\"ahler metrics  \ on $V$, the first one is smooth
with potential $\f^1=|u|^2+ |v|^2+ |w|^2$, induced by the euclidean K\"ahler metric of  $\C^3$, the second one
is the  K\"ahler current whose potential is $\f^2= |u|+|v|$. On $V^{reg}$ it is the quotient
of the euclidean metric restricted to $\C^2-\{0\}$.    Near $0$, $dd^c\f^2 \gg dd^c \f^1$.
\end{exa}

The metric $dd^c\f^2$ is an example of an orbifold K\"ahler metric on $V$. 
The results of [Y] extend without
major modifications to K\"ahler orbifolds. For instance, in each K\"ahler class of a nodal K3
surface there is a unique Ricci flat orbifold metric.

\subsubsection*{Chern-Weil forms and hermitian metrics}

Let $\mathcal{PH}_V$ be the sheaf of real-valued pluriharmonic functions on $V$. 
By definition, a closed (1,1)-form on $V$
is a  section of the sheaf ${\mathcal C}^{\infty}_V/ \mathcal{PH}_V$. We have the exact sequence:
$$
{\mathcal C}^{\infty}(V)\to \Gamma(V, {\mathcal C}^{\infty}_V/ \mathcal{PH}_V) 
\buildrel{[ \ \ .\ \ ]}\over{\longrightarrow} H^1(V,\mathcal{PH}_V) \to 0. 
$$

A class in  $H^1(X,\mathcal{PH}_X)$ will be called  K\"ahler, if it is 
in the $[\ \ . \ \ ]$ image of a smooth K\"ahler metric.

\begin{rqe} Assume $X$ is smooth.
A class $[\omega]$ in $H^1(X,\mathcal{PH}_X)$ will be called numerically base point free iff
there exists  a proper surjective holomorphic mapping $X\to Y$, $Y$ normal, such that $[\omega]$ is the
pull back of a K\"ahler class on $Y$.
This is in principle a stronger condition
than being cohomologous to a smooth semipositive closed $(1,1)$-form, although no counterexample seems to be known.

If $X$ is projective and $[\omega]\in NS_ {\mathbb R}(X)$,  $[\omega]$ is numerically base point free class iff it is semiample.

In the non-big case (i.e.: $\int_X \omega^n =0$), it is straightforward to construct 
semi-K\"ahler classes that
are not numerically base point free (e.g. on complex tori).
\end{rqe}

Let  $L$ be a holomorphic line bundle on $V$.
The notion of smooth hermitian metric on
 $(V,L)$ is defined as in the smooth case.
Let $h$ be such a metric on $(V,L)$.

Let $s\in H^0(U, L)$ be a nowhere zero local holomorphic section of $L$ (a local generator of $L$) 
defined over the open subset $U\subset V$. Set  $e^{-\f_ s}:= ||s||_h^2$,
where $\f_s$ is a ${\mathcal C}^{\infty}$-smooth function on $U$. 
The current $ dd^c\f_s$ is a smooth closed (1,1)-form on $V$
which does not depend on $s$; it is
a semi-K\"ahler current if $\f_s$ is psh.

More generally,  let $(U_i)_i$ be an open covering of
$V$ and $\sigma_i\in H^0(U_i, \Ox(L))$ a local generator of $L$.
Let $\f_i=\f_{\sigma_i}$. The datum $(U_i,\f_i)$ defines  a smooth closed (1,1)-form on $V$.

\begin{defi}
The  Chern-Weil form of $(V,L,h)$ (or of $h$) is the $\sim$-equivalence class of  the data $(U_i,\f_i)$
constructed above. We will denote it by $c_1(L,h)$.
\end{defi}

It is immediate that $[c_1(L,h)]$ is independent of   $h$. 
Hence there is a linear map $c_1: Pic(V)\to H^1(V,\mathcal{PH}_V)$.
The connection with the more widely known smooth case is made by the observation that,
if $X$ is a compact K\"ahler manifold, $H^{1,1}(X,\R)= H^1(X,\mathcal{PH}_X)$.

\begin{pro}
Let  $V$ a compact normal complex analytic  variety.
 
The space $H^1(V,\mathcal{PH}_V) $ is finite dimensional.

Let $L$ a holomorphic line bundle on $V$. Every representative of $c_1(L)$ in $H^1(V,\mathcal{PH}_V)$
is the Chern-Weil form of a smooth  hermitian on $L$.

 If there exists a smooth hermitian metric $h$
such that $c_1(L,h)$ is K\"ahler, then $V$ is projective-algebraic and $L$ is  ample.
\end{pro}

\begin{proof}
The most difficult task is to show  that, in the last assertion, $V$ is Moishezon. This
 follows from Siu's solution
of the  Grauert-Riemenschneider conjecture [Siu].
\end{proof}
\vskip.1cm

 A singular metric on  $L$ is an expression $h= e^{-\f} h_{sm}$, $\f$ being a locally smooth + psh function and  $h_{sm}$ a smooth   hermitian metric.
 Its  Chern-Weil form is the  quasi-positive current $c_1(L,h_{sm})+dd^c\f$.

\section{Adapted volume forms}

\subsection{Monge-Amp\`ere equations on normal K\"ahler spaces}

Let $\Omega$ be a smooth K\"ahler metric on $V$.  A classical result of P. Lelong states that if
$U$ is relatively compact in $V$,  then  $U^{reg}$
is of finite volume with respect to the smooth
volume form $\Omega_{reg}^n$. 

This has been generalized by E.Bedford and A.Taylor in [BT], where
the authors study Monge-Amp\`ere
measures for locally bounded psh functions. Since these measures do 
not charge proper analytic subsets, we obtain:

\begin{pro}
Let $\Omega$ be a semi-K\"ahler current with  $L_{loc}^{\infty}$ potentials on  $V$.
The Monge-Amp\`ere measure $\Omega_{reg}^n$ is well defined on $V_{reg}$ and satisfies
$ \int_{U^{reg}} \Omega_{reg}^n< \infty, 
\text{ for all relatively compact subset } U \subset V.$

For any resolution $\pi: X\to V$, the Monge-Amp\`ere measure $(\pi^*\Omega)^n$
is well defined on $X$ and  satisfies $\pi_*(\pi^*\Omega)^n = j_*\Omega_{reg}^n$.
Moreover if $\bar \pi: \bar X\to V$  is a  resolution  dominating $\pi$
(i.e. $\bar \pi= \pi \circ \psi$
for some bimeromorphic proper holomorphic map
$\psi: \bar  X' \to X'$), then $\psi_* (\bar \pi^* \Omega )^n= (\pi^*\Omega)^n$.
\end{pro}

The measure $\pi_*(\pi^*\Omega)^n$ is thus well defined on $V$ and independent of the
choice of resolution.
We will call it the {\bf Monge-Amp\`ere measure of
$\Omega$} and denote it  by $\Omega ^n$.
The mass of this measure only depends on the cohomology class of $\Omega$,
as follows again from [BT]:

\begin{lem} Assume $V$ is compact.
 Let $\Omega_1,\Omega_2$ two semi K\"ahler currents with   $L_{loc}^{\infty}$ potentiel on  $V$.
If they are cohomologous, i.e. $\Omega_1=\Omega_2+dd^c \f$ for some $\f\in L^{\infty}(X)$,
then $\int _V \Omega_1^n= \int _V \Omega_2^n$.
\end{lem}

We can now reformulate some of our previous results. 

\begin{thm}\label{sma} 
Let $V$ be a $n$-dimensional compact normal K\"ahler space and $\Omega$ be a smooth 
K\"ahler form on $V$. 
Then for every $f\in L^p(V,\Omega^n)$,  $p>1$, such that 
$\int_V f \Omega^n= \int_X \Omega^n$, there is a unique $\f\in L^{\infty}(V)$ such that 
$$
(\Omega+dd^c\f)^n = f\Omega^n \text{ and } \sup_V \f =-1.
$$ 
\end{thm}

\begin{proof}
Let $\pi: X\to V$ be a resolution of $V$. We may define a semipositive big smooth form
on $X$ by $\omega=\pi^* \Omega$. By Theorem \ref{mas} and Proposition \ref{lp} we can solve
uniquely
$(\omega+ dd^c \bar \f)^n = f \circ \pi \omega^n$ where $\bar \f$ is a continuous function on $X$
such that $\omega+dd^c\bar \f$ is semipositive. Let $F$ be a fiber of $\pi$ and $i: F\to X$
the inclusion map. 
$F$ is connected by Zariski's main theorem. Furthermore $i^*\omega+dd^c i^*\bar\f$
is semipositive on $F$. Since $i^*\omega=0$, it follows that $i^*\f$ is a psh function 
on $F$. Hence $i^*\bar\f$ is constant. This implies that $\bar\f = \f\circ \pi$
where $\f$ is a bounded u.s.c. function on $V$. We do have $(\Omega+dd^c\f)^n = f\Omega^n$.
\end{proof}

\subsection{Adapted measures on log terminal K\"ahler spaces}

Let $V$ be a $n$-dimensional Gorenstein K\"ahler space and $\Omega$ be a smooth 
K\"ahler form on $V$. Fix  $x\in V$ and let $\alpha$ be a local generator of 
$\omega_V$ defined over an open subset $x\in U$; then
$v=c_n \alpha \wedge \bar\alpha$ is a 
positive definite volume form on $U^{reg}$, for an appropriate choice 
of the constant $c_n=\sqrt{-1}^{n} (-1)^{\frac{n(n+1)}{2}}$.
 
When $V$ is merely $\Q$-Gorenstein of finite index $N$, 
we choose $\beta$ a local generator of
$\omega_V^{[N]}$ defined over an open subset $x\in U$ and we set
$$
v=v_{\beta}=\left( \sqrt{-1}^{Nn} (-1)^{N\frac{n(n+1)}{2}}\beta \wedge \bar\beta \right)^{\frac{1}{N}}.
$$
This is a positive definite volume form on $U^{reg}$. 

Our next observation is that log terminal singularities are the worst singularities
we can allow in order to {\it globally} solve Monge-Amp\`ere equations associated
to volume forms on $V$.

\begin{lem}\label{lpd}
For every $U_1\subset\subset U$, $\int_{U_1^{reg}} v<\infty$ iff $X$ is log terminal.
 
 If $V$ is log terminal, then the Radon measure $\mu= j_*v$ satisfies 
$\mu= f\Omega^n$ with $f\in L^{1+\e}(U_1,\Omega^n)$ for some $\e>0$.
\end{lem}

\begin{proof}
Let $\pi :X\to V$ be a log resolution. 
Write $K_{X}\cong \pi^* K_V + \displaystyle\sum a_E E$. 
Since $exc(\pi)$ has simple normal crossings, at every $P\in E=exc(\pi)$ 
there are local coordinates $(z^i)_{i=1,...,n}$ such that $E$ is described by the equation
 $z^1\ldots z^q=0$. Let $E_j$ be the divisor $z_j=0$. We have:
 $\pi^*v= \prod_{j=1}^q |z^j|^{ 2a_{E_j}} d\lambda$ where $d\lambda$ 
 is a Lebesgue measure on $X$, hence
the measure $\pi^*v$ has finite mass near $P$ iff $\forall j, a_{E_j} >-1$. 
Thus $\int_{U_1^{reg}} v<\infty$ iff $\forall E, \ a_E > -1$. 

Let $f_1$ be the density of $\pi^*v$ with respect to $d\lambda$.  
Since $f_1$ is comparable to $\prod_{j=1}^q |z^j|^{ 2a_{E_j}}$ near $P$,
it follows that $f_1$ belongs actually to $L^p(X,d\l)$ for some $p>1$ 
when $X$ is log terminal.
 
Let $D=1/f$ be the density of $\Omega^n$ with respect to $v$.
We will see here below that $D$ is bounded but it might have zeroes on $E$,
hence $f$ is unbounded in general. However we will show that
$f \circ \pi \in L^{\a}(X,d\l)$ for $\a>0$ small enough, hence it
follows from H\"older's inequality (as in the proof of lemma 3.2) that 
$$
\int_{U_1^{reg}} f^{1+\e}\Omega^n= 
\int_{ \pi^{-1}U_1^{reg}} f^{\e}  f_1  d\lambda <+\infty
$$
if $\e>0$ is small enough.

Fix $x\in V$ and let $i:U_x\to \C^m$ be a local embedding of a neighborhood 
$U_x$ of $x$. We consider the 
 $\left( 
 \begin{array}{c}
  m \\
   n
    \end{array}
     \right)$
 $n$-forms
on $U_x^{reg}$ $du^I=du^{i_1}\wedge \ldots du^{i_n}$, 
where $(u^i)$ is a set of affine coordinates on $\C^m$. 
Observe that $\Omega^n$ is comparable to $\sum_I v_{du^I}$ 
\footnote{Note that the formula for $v_{\beta}$ makes sense 
even if $\beta$ is not a local generator.}. 
Since $\beta$ is a local generator at $x$ of $\omega_V^{[N]}$, we
have $(du^I)^N = f_I \beta$ where $f_I\in \mathcal{O}_{V,x}$  is the germ
of a holomorphic function at $x$.
Therefore  $\Omega^n$ is comparable to $\sum_I |f_I|^{\frac{2}{N}}v$,
hence $D$ is comparable to $ [ \sum_I |f_I|^{\frac{2}{N}}]^{-1}$
near $x$.
 
 The functions $(f_I)$ generate an ideal $\mathcal I_x\subset \mathcal{O}_{V,x}$. 
Actually, the construction can be globalized
to provide a coherent ideal sheaf $\mathcal I \subset \mathcal{O}_V$ 
cosupported on $V^{sing}$. 
 
 We may assume [Hi], [BM] that $\pi: X\to V$ is a log resolution of $(V,\mathcal I)$, 
namely a log resolution 
 of $V$ with the additional property that the ideal sheaf $\pi^{-1}\mathcal{I} . \mathcal{O}_X$
 which is the ideal sheaf of $\mathcal{O}_X$ generated by the family of holomorphic functions
 $(\pi^* f_I)_I$, satisfies 
$\pi^{-1}\mathcal{I} . \mathcal{O}_X = \mathcal{O}_X (-\sum Nb_E E)\subset \mathcal{O}_X$
 where $N.b_E\in \N$ is a positive multiplicity attached to any exceptional divisor of $\pi$. 
 
 In local coordinates near $P\in X$, 
$\pi^*D$ is comparable to $\prod_j |z_j| ^{ 2b_{E_j}}$, hence
$
\pi^*(f_1D^{-\e}) \text{ is comparable to } \prod_j |z_j|^{2(a_{E_j} -\e b_{E_j})}.
$
It follows that for every relatively compact subset 
$U_1\subset\subset U, \ \ f\in L^{1+\e}(U_1,\Omega^n)$ 
iff $\forall E, \ \  \pi(E)\cap  U \not=\emptyset  \Rightarrow a_E- \e b_E >-1$. 
\end{proof}

\begin{defi}
 Assume $V$ has only log terminal singularities.  
A positive definite adapted measure on $V$ is a positive Radon measure
  locally of the form $e^f. v$
 where $f$ is a bounded measurable function. 
 A positive definite adapted measure has ${\mathcal C}^0$, 
${\mathcal C}^{\alpha}$, ${\mathcal C}^{\infty}$ density if so is $f$.
\end{defi}
   
\begin{rqe}
It follows from lemma \ref{lpd} that
if $V$ is $\Q$-Gorenstein but has non log terminal singularities, $v$ is a volume form on 
$V^{reg}$ but does not extend to a measure on $V$.
\end{rqe}

\subsection{Adapted volume forms for klt K\"ahler pairs}

We will be briefer since pairs are mainly of interest to MMP practitioners.
The key definition for us will be: 

\begin{defi} A pair $(V,\Delta)$ is klt iff $K_V+\Delta$ is $\Q$-Cartier 
and if for any log-resolution 
$\pi:X\to V$
of $(V,\Delta)$,  we have the numerical equivalence of Cartier divisors:
$$
N(K_X+\Delta')\cong \pi^*N(K_V+\Delta)+ \sum_{E \  exc.} Na_E E
$$
with $a_E>-1$, $\Delta'$ the proper transform of $\Delta$ in $X$ (same multiplicities)
and $N$ is an integer such that $N(K_V+\Delta)$ is Cartier. 
\end{defi}

Thus a variety $V$ has only klt singularities iff $(X,\emptyset)$ is klt. 
 
Let $\beta$ be a local generator of $\Ox(N(K_V+\Delta))$. 
Then $\beta_{V^{reg}}$ can be viewed as a meromorphic N-canonical form
with a pole of order $Nd_i$ on $E_i$ where $\Delta=\sum_i d_i E_i$
is the decomposition of $\Delta$ into prime divisors.
Thus $v_{\beta_{V^{reg}}}$ defines a volume form with poles on $V^{reg}$,
namely
$v_{\beta_{V^{reg}}}$ is comparable to $ \prod_i |\sigma_i|^{-2d_i} d\l$,
where $\sigma_i$ denotes the canonical section of $\mathcal{O}(E_i)$. If $v_{\beta_{V^{reg}}}$
is a finite measure then $d_i<1$, but the converse is not true. We have the following 
staightforward extension of lemma \ref{lpd}:

\begin{lem} Let $j': V-\cup_i E_i \to V$ be the canonical inclusion.
$j'_* v_{\beta}$ is a well defined  Radon measure on $V$ iff $(X,\Delta)$ is klt. 
\end{lem}

The definition of an adapted measure for a klt pair is left to the reader.

\section{Singular K\"ahler-Einstein metrics}

 \subsection{Singular Ricci curvature}

\subsubsection*{The smooth case}

The link between Monge-Amp\`ere equations and K\"ahler-Einstein metrics is provided
by the following classical

 \begin{lem} \label{ric}
 Let $X$ be complex manifold, let $h$ be a smooth hermitian metric
on $\omega_X$ and $\Omega$ a K\"ahler form such that
 $\Omega^n= v(h)$. The Ricci curvature divided by $2\pi$ of $\Omega$ is the Chern-Weil form 
 $-c_1(K_X,h)$. 
 \end{lem}
 
 \subsubsection*{Adapted measures and hermitian metrics on the canonical sheaf}

Assume $V$ is compact with only canonical singularities, has index  $N$ and let $h^N$ be a smooth hermitian metric
 on $\omega_V^{[N]}$. Let $\beta$ be a local generator of  $\omega_V^{[N]}$. 
 Define $v_{\beta}(h)$ to be the volume form on $V^{reg}$: 

$$ 
v_{\beta}(h) =\left( \sqrt{-1}^{Nn} (-1)^{N\frac{n(n+1)}{2}}
\frac {\beta \wedge \bar\beta}{\|\beta \|^2_{h^N}}
 \right)^{\frac{1}{N}}
$$ 

Since $v_{\beta}(h)$ is independent of $\beta$,  this expression defines an
 adapted measure $v(h)$ with ${\mathcal C}^{\infty}$ density on $V$. 
 
 Now, let $h^N_{sing}=e^{-N\chi} h^N$ be a singular metric on $\omega_V^{[N]}$.
The Chern-Weil form  $c_1(\omega_X^N,h_{sing}^N)$ 
is then well defined as a quasipositive current. Since $h_{sing}^N$ has locally ${\mathcal C}^{\infty}$+psh potentials
$\chi$ is locally bounded above
and the above formula defines a measure $v(h_{sing})=e^{\chi} v(h)$ on $V$ such that   $\frac{v(h_{sing})}{ v(h)}\in L^{\infty}_{loc}$. 
In particular 
$$
\frac{v(h_{sing}) }{\Omega^n} \in L^{1+\e} (V,\Omega^n)
\text{ for } \e>0 \text{ small enough}. 
$$
We have $c_1(K_X,h_{sing})= c_1(K_X,h)+ dd^c\chi$,
where $c_1(K_X, h):= \frac{1}{N} c_1(\omega^N_X, h^N)$.

   \begin{defi}
 Assume $V$ has only canonical singularities.  An adapted measure on $V$ is a positive Radon measure
  locally of the form $e^f. v$
 where $f$ is locally given as the sum of a psh and a smooth function on $V$. 
 An adapted measure has ${\mathcal C}^0$, resp. ${\mathcal C}^{\alpha}$, resp. ${\mathcal C}^{\infty}$ density if $e^f$ is ${\mathcal C}^0$, resp.  ${\mathcal C}^{\alpha}$, resp. ${\mathcal C}^{\infty}$.
   \end{defi}
   
   The definition has the virtue of generalizing the usual  equivalence between smooth 
metrics on the canonical sheaf of a manifold and
positive definite volume forms to singular metrics and  log terminal spaces.
This suggests the following:

\begin{defi}\label{ske}
Let $V$ be a $\Q$-Gorenstein K\"ahler normal $n$-dimensional complex space with only 
log terminal singularities.
Let $\Omega$ be a semi-K\"ahler current with $L^{\infty}_{loc}$ potential 
and adapted Monge-Amp\`ere measure. 
Let $h$ be the singular metric on the canonical sheaf such that $\Omega^n = v(h)$. 
We define 
$$
Ric(\Omega):=-c_1(K_V,h),  
$$
where the equality is to be taken in the sense of currents.

$\Omega$ will be called a singular K\"ahler-Einstein metric
if $Ric(\Omega)=c \Omega$ for some $c \in \R$.
\end{defi}

\subsection{Singular Ricci flat metrics}

\begin{defi}
Let $V$ be a K\"ahler space with only canonical singularities.  $V$ is said to be $\Q$-CY, iff 
there is some multiple $N'$ of $index (X)$  such that $H^0(V,\omega^{[N']}_V)=\C \alpha$,
where $\alpha$ is
a global generator of $\omega^{[N']}_V$. 
\end{defi}

\begin{thm} 
Assume $V$ is a compact $\Q$-CY K\"ahler space. 
Let $\Omega$ be a smooth K\"ahler metric on $V$. 
Then there is a unique semi-K\"ahler current with locally bounded potential 
and adapted Monge-Amp\`ere measure 
$\Omega'=\Omega +dd^c\f$, such that 
$$
(\Omega+dd^c\f)^n = C v_{\alpha}
\text{ and } \sup_V \f =-1,
$$
where   
$\int_V \Omega^n=C \int_V (-1)^n v_{\alpha}$.
 
Furthermore, if $V$ is projective-algebraic and $[\Omega] \in NS_{\R}(V)$, then $\Omega+ dd^c\f$ is smooth on $V^{reg}$ where it defines a bona fide 
Ricci flat metric. 
\end{thm}

\begin{cor}
In each cohomology class of a smooth K\"ahler form,
there is a unique singular Ricci flat metric.
\end{cor}

\begin{proof}
This follows straighforwardly from Theorems \ref{sma}, 3.6, Lemma \ref{lpd} and Definition \ref{ske}.  
\end{proof}
\vskip.1cm

 \begin{exa}
 A nodal quintic threefold is $\Q$-CY and has not quotient singularities, so the orbifold
 method of [Ko] does not work. 
 \end{exa}

\subsection{Singular K\"ahler-Einstein metrics of negative curvature}

\begin{thm} \label{tg}
Let $V$ be a general type projective algebraic variety with only canonical singularities
such that $K_V$ is ample.
 Let $h^N$ be a smooth hermitian metric on  $\omega^N_V$ such that $\Omega= c_1( K_V, h)$
 is a smooth K\"ahler form on $V$.
 
 There is a unique  $\f \in L^{\infty}(V,\R)$ such that:
 \begin{enumerate}
 \item $\f$ is $\Omega$-psh.
 \item $\Omega+ dd^c\f$ semi K\"ahler current with  $L^{\infty}$ potential.
 \item $(\Omega+dd^c\f)^n=e^{\f}v(h)$. 
  \end{enumerate}
Consequently $\Omega+dd^c \f$ is the unique singular KE metric on $V$ of negative curvature in
  the canonical class of $V$. 
The current $\Omega+dd^c \f$ has locally bounded potentials
and is smooth on $V^{reg}$ where it defines a bona fide 
KE metric. 
  \end{thm}

\begin{proof}
This is a consequence of Theorems \ref{mma}, 4.4, and Definition \ref{ske}. 
\end{proof}

 \begin{rqe}
 Thanks to Theorem \ref{Reid}, for
 $X$ a projective algebraic manifold of general type such that 
$R( X):=\oplus_{n\in\N} H^0( X, \mathcal{O}_X (nK_{ X}))$
 is finitely generated, $X$ has a unique birational model $V$ such that the above hypotheses hold. 
Thus we have a 
 birational map $\pi: X\dashrightarrow V$ which
is well defined outside an indeterminacy locus $S$ of codimension
 $\le 2$.   In particular $\pi^* (\omega+dd^c\f)$ is a closed positive current on $X-S$ 
that extends to a 
 closed positive current $T$ on $X$ itself. 
The current $T$ defines a KE metric on $X-S$. It needs not be a singular 
 KE metric on $X$ though, since its potentials may have logarithmic poles on $S$, in fact algebraic singularities of the form $\alpha \log (\sum |f_i|^2)+ O(1)$
 $f_i$ holomorphic and $\alpha\in \Q_{>0}$. Moreover, $T$  lies in the canonical class
 of $X$ iff $X$ is a smooth minimal model as in [Ts].
 \end{rqe}

\subsubsection*{Connection with [Ts]}
 
Let $X$ be a complex projective manifold such that $K_X$ is nef and big. 
Let $\Omega$ be a smooth K\"ahler metric on $X$ and consider the 
K\"ahler-Ricci flow
$$ 
\frac{\partial \Omega_t}{\partial t} = -Ric(\Omega_t) -\Omega_t, \ \ \Omega_0=\Omega.
$$

In [Ts], it was proved that this flow has a global solution for all time $t\in [0,\infty[$, and 
an argument was given, recently  fully completed in [TZ], 
to the effect that $\Omega_t$ converges to a closed positive current $T_{KE}$, 
independent of $\Omega$, which defines a smooth 
K\"ahler-Einstein metric outside the exceptional divisor $E$ of the holomorphic 
bimeromorphic map $X\to X_{can}$. Its 
potential  satisfies  the 
Monge Amp\`ere \'equation considered in Theorem \ref{tg} outside $E$. 
It follows from proposition 4.4 that the current $T_{KE}$ coincides with
the solution produced by Theorem 7.8.

The independant work  [TZ] gives a proof of the following 
 properties, already conjectured by [Ts], that $T_{KE}$ has locally
 bounded potential and satisfies the degenerate Monge-Amp\`ere
 equation considered in Theorem \ref{tg}.

 \begin{exa}
 A nodal sextic threefold is of general type, Gorenstein, terminal, is its own canonical model, 
has no smooth minimal model and 
does not have quotient singularities.  Therefore  the orbifold
 method of [Ko] does not work and  [Ts]
 does not apply. 
 \end{exa}

\subsection{Singular KE metrics on klt pairs}
 
Let us now state the immediate generalization to klt pairs.

 \begin{defi}
Let $(V,\Delta)$ be a klt compact K\"ahler pair. 
 
The pair $(V,\Delta)$ is said to be $\Q$-CY, iff 
there is some multiple $N'$ of $index (X,\Delta)$  such that 
$H^0(V,\Ox (N'(K_V+\Delta)))=\C \alpha$  where $\alpha$ is
a global generator of $\Ox (N'(K_V+\Delta))$. 
 
The pair $(V,\Delta)$ is canonically polarized
iff $K_V+\Delta$ is ample. 
\end{defi}

\begin{thm}
Let $(V,\Delta)$ be a klt compact K\"ahler pair. 

If $(V,\Delta)$ is $\Q$-CY it carries a 
singular Ricci flat metric with adapted volume form in any K\"ahler class of $V$, 
this current being smooth outside $\Delta\cup V^{sing}$ if $V$ projective and the K\"ahler 
class is rational. 

If it is canonically polarized it
carries a unique singular KE metric in the cohomology class of $K_V+\Delta$, regular outside
$\Delta\cup V^{sing}$. 

Furthermore, let $V^o$ be the largest open subset of $V^{reg}$ such that $\Delta \cap {V^o}$ has snc support and 
multiplicities of the form $1-\frac{1}{n}$ with $n\in \mathbb N^*$. Then, the singular KE metric becomes smooth on the
stack $[V^o, \Delta \cap {V^o}]$. 
\end{thm}

\begin{proof}
For regularity on the smooth locus, we need the full statement of Theorems \ref{reg} and \ref{reg2}, poles included.
\end{proof} 

\vskip.3cm

\begin{ack} We would  like to thank Z. Blocki, S. Boucksom, 
A. Chiodo, J.P. Demailly, J. Keller, S. Kolodziej, M. Paun and B. To\"en for useful conversations
and C. Simpson and Y.T. Siu for  inspiring remarks. We are grateful to the referee for his penetrating remarks.
\end{ack}

\vskip .3cm

Philippe Eyssidieux

Institut Fourier - UMR5582

100 rue des Maths, BP 74

38402 St Martin d'Heres (FRANCE)

eyssi@fourier.ujf-grenoble.fr

\vskip.2cm

Vincent Guedj 

LATP

UMR 6632, CMI, Universit\'e de Provence

39 Rue Joliot-Curie

13453 Marseille cedex 13 (FRANCE)

guedj@cmi.univ-mrs.fr

\vskip.2cm

Ahmed Zeriahi

Institut de Math\'ematiques de Toulouse (IMT)

Universit\'e Paul Sabatier, 118 route de Narbonne

31062 TOULOUSE Cedex 04 (FRANCE)

zeriahi@math.univ-toulouse.fr

\end{document}